\documentclass[12pt]{article}
\usepackage{amsfonts,amssymb,amsmath,amsthm}

\usepackage{xcolor}
\definecolor{darkgreen}{rgb}{0.,0.3,0.}

\usepackage{graphicx}

\newcommand{\R}{{\mathbb R}}

\newcommand{\Z}{{\bf Z}}

\newcommand{\D}{{\bf D}}

\newcommand{\Fix}{{\rm Fix}}
\newcommand{\be}{{\bf e}}
\newcommand{\ou}{{\rm out}}
\newcommand{\ini}{{\rm in}}
\newcommand{\bx}{{\bf x}}
\newcommand{\bg}{{\bf g}}
\newcommand{\bh}{{\bf h}}
\newcommand{\bv}{{\bf v}}

\newcommand{\bu}{{\bf u}}

\newcommand{\rd}{{\rm d}}
\newcommand{\cB}{{\cal B}}

\newcommand{\qtilw}{{\widetilde{PQ}}_w}

\newcommand{\qw}{PQ_w}
\newcommand{\qz}{PQ_z}
\newcommand{\rw}{R_w}
\newcommand{\rz}{R_z}

\newcommand{\dpt}{\displaystyle}

\newtheorem{theorem}{Theorem}[section]
\newtheorem{lemma}[theorem]{Lemma}
\newtheorem{Cy}[theorem]{Corollary}
\newtheorem{definition}[theorem]{Definition}
\newtheorem{remark}[theorem]{Remark}

\begin{document}
\begin{center}
\mbox{}
\bigskip

\noindent{\large \bf Stability of a heteroclinic network and its cycles:\\
a case study from Boussinesq convection}\\
\bigskip

{\small \bf Olga Podvigina}\\
{\footnotesize Institute of Earthquake Prediction Theory\\
and Mathematical Geophysics,\\
84/32 Profsoyuznaya St, 117997 Moscow, Russian Federation \\
email: olgap@mitp.ru}

\bigskip
{\small \bf Sofia B.S.D. Castro}\\
{\footnotesize Centro de Matem\'atica and Faculdade de Economia\\
Universidade do Porto\\
Rua Dr. Roberto Frias, 4200-464 Porto, Portugal\\
email: sdcastro@fep.up.pt}

\bigskip

{\small \bf Isabel S. Labouriau }\\
{\footnotesize Centro de Matem\'atica da Universidade do Porto\\
Rua do Campo Alegre 687, 4169-007 Porto, Portugal\\
email: islabour@fc.up.pt}

\end{center}

\begin{abstract}
This article is concerned with three heteroclinic cycles forming a  heteroclinic network in $\R^6$.
The stability of the cycles and of the network are studied.
The cycles are of a type that has not been studied before, and provide an illustration for the difficulties arising in dealing with cycles and networks in high dimension.
In order to obtain information on the stability for the present network and cycles, in addition to the information on eigenvalues and transition matrices, it is necessary to perform a detailed geometric analysis of return maps.
Some general results and tools for this type of analysis are also developed here.
\end{abstract}

\section{Introduction}

In this article we derive stability conditions for a specific heteroclinic network in $\R^6$, as well as for its cycles. This network is of a type that has not been studied before and has features that distinguish it clearly from what is discussed in the literature.
This case study  both provides a starting point for further general stability results and illustrates the difficulties arising in the study of higher-dimensional more general networks.

Recall that the smallest dimension where a robust heteroclinic cycle can exist
is $n=3$. Robust heteroclinic cycles existing in $\R^3$ have been known for a long
time, going as far back as the work of dos Reis \cite{dosReis} and Guckenheimer and Holmes \cite{GuckenheimerHolmes}; 
the list of possible cycles is short.
In $\R^4$
the situation becomes more complex.
However, general results on heteroclinic cycles and networks in $\R^4$ are
known in the literature, starting with that by Krupa and Melbourne \cite{KrupaMelbourne1}.
In \cite{KrupaMelbourne1} the term ``simple'' was
attributed to robust heteroclinic cycles emerging in $\Gamma$-equivariant
systems in $\R^4$, such that, in particular, heteroclinic connections belong
to planes that are fixed point subspaces for subgroups of $\Gamma$.
Depending on how the subgroups act on $\R^4$, simple heteroclinic cycles
were further subdivided into types A, B and C.
The definitions of simple and type A cycles were extended to higher
dimensions in\cite{KrupaMelbourne2}, \cite{op12} also in terms of how the subgroups act on certain
invariant subspaces, while in this spirit the cycles of types B and C were
generalised as type Z in \cite{op12}.

In $\R^5$, the list of finite subgroups of O(5) is known: it is a union of finite
subgroups of O(4) and a few other subgroups \cite[ArXiv version]{Mochia}, therefore
it is likely that heteroclinic cycles existing in $\R^5$ are not very different from the
ones in $\R^4$.
This is certainly the case for homoclinic cycles \cite{op13,Sottocornola}.
Some instances of heteroclinic cycles in $\R^6$ were considered in  the literature \cite{AguiarCastro,GdSC2017},
however no general results are yet available.  Systematic ways of constructing, not necessarily simple, heteroclinic cycles in any dimension have been established in \cite{AshwinPostlethwaite,Field}.

Concerning their stability properties, heteroclinic cycles in $\R^3$ are either asymptotically stable or
completely unstable and the conditions for asymptotic stability are
trivial. In $\R^4$, cycles that are not asymptotically stable can be stable in
a weaker sense, namely essentially \cite{Melbourne1991} or fragmentarily asymptotically stable
\cite{op12}. Of the two, essential asymptotic stability is the strongest.
In \cite{pa11,op12,pc16}, conditions for stability  for simple and pseudo-simple cycles in $\R^4$ are obtained 
 from the eigenvalues of the Jacobian at the nodes of the cycle and/or
from eigenvalues and eigenvectors of so-called transition matrices.
For cycles that are not simple but for which the
transitions  along connections behave as permutations, analogous tools can be
used to establish stability properties \cite{GdSC2017}.
 The stability of heteroclinic cycles may also be studied by
making use of Lyapunov functions, as in \cite{HofbauerSigmund} in the context of
population dynamics (non-simple cycles).
The network in the present case study is not simple and is different from those
considered in \cite{GdSC2017} ,calling for different techniques in the study of stability.

Loss of stability, as well as stability itself, is the starting point for further studying the
dynamics near the heteroclinic cycle or network and has been pursued by
several authors. A selection of examples is given in
\cite{Lohse15,Postlethwaite,PostlethwaiteDawes}. This further development
is out of the scope of the present article.

Before addressing the case study, we prove generic results that apply to any robust heteroclinic network in an Euclidean space of any finite dimension. The main general result is on the (lack of) asymptotic stability of  networks consisting of a finite number of one dimensional  connections.

The network in the case study is such that neither eigenvalues of the Jacobian at its nodes nor transition matrices provide complete information about stability.
To overcome this, we obtain stability results for the fixed points of several families of maps
that
have the generic analytic form of simplified return maps to cross-sections to connections in a cycle or network.
These results may be useful in the study of generic robust heteroclinic cycles or networks.
The stability results we establish for the fixed points of these maps are crucial for the study of the stability of our particular network.

 The network in the case study has been described in \cite{clp} in the context of a convection problem.
We obtain fragmentary asymptotic stability conditions for this network and for its cycles in the following four steps:
\begin{itemize}
\item[(a)]
obtain a first return map $g$ as the composition of local maps around nodes and global transition maps;
\item [(b)]
obtain from $g$ a reduced map $h$, defined in a lower dimension;
\item[(c)]
find stability conditions for fixed points of $h$;
\item [(d)]
show that  the stability conditions for $h$ coincide with the stability conditions for $g$.
\end{itemize}
Then we obtain more information:
\begin{itemize}
\item[(e)]
deriving conditions for  essential asymptotic stability from stability indices.
\end{itemize}

Step (a) is algorithmic and well known, although it may yield cumbersome expressions when either the phase space dimension or the length of the cycle is large.
The other steps are non-standard.
Our study indicates that they may always be done in roughly the same way, but with a procedure that has to be reinvented for each case.

Step (b) is not easy but maybe a general formulation is possible, although complicated.

Step (c) is certainly very difficult and we have no hope of generalising it, in particular, for lack of a general form for $h$.
We make a geometric analysis of the stability, adapting to each case the  results on the stability of fixed points of general maps.

Step (d) perhaps can be given a general proof, but certainly it will be highly non-trivial and not worthwhile trying since one does not have a generalisation for step (c).

Step (e) is the only one that is not so difficult in our case, once the others were done.
It is not clear what would happen in other cycles or networks but addressing a more general case is beyond the scope of this article.
\bigbreak

In the cases of type A or Z  cycles, the stability can be decided  from information on eigenvalues and eigenvectors of the linearisation at nodes and of  transition matrices.
This then can be used to obtain general results for these types.
For  other cycles in $\R^n$, in particular for larger $n$, the linear information has to be used in a more involved way.
Steps (a) to (e) above provide a heuristic approach for deciding the stability of a heteroclinic object  in $\R^n$, in the cases where knowing the eigenvalues and eigenvectors is not sufficient to decide stability.
Our example leaves little hope of finding general conditions for stability that may be stated in a simple way, except for very specific classes of cycles.

We finish this section with a short description of the network and its stability. In the next section we provide some technical background. Section~\ref{stability_results} provides generic stability results, while Subsection~\ref{subDescription} describes the network which is our main concern,
 with details in  Appendix~\ref{AppTable}.
In the remainder of Section~\ref{CaseStudy} we address the stability of individual cycles and  of the network as a whole.
 The final section concludes.
 
 We consider a network that is a union of three heteroclinic cycles.
The network emerges in a twelve-dimensional dynamical system obtained by
the center manifold reduction from the equations of plane layer
Boussinesq convection with a hexagonal periodicity lattice, see \cite{clp}.
The symmetries, and hidden symmetries,
of the problem allow for a further reduction to six dimensions.
The symmetry group of the system implies the existence of several
 one- and two-dimensional flow-invariant subspaces, that were identified
in \cite{clp} where we derived conditions for existence of
structurally stable heteroclinic cycles.
These were presented as
inequalities involving the normal form coefficients.
Numerical simulations in \cite{clp} illustrate the behaviour of trajectories
near two of these cycles and indicate that the third cycle is
completely unstable, at least for the considered values of coefficients.
In the present paper our sole assumption for the study of stability is the
standard one that the equilibria
involved in the network are stable in the transverse directions, i.e. all
eigenvalues not related to outgoing heteroclinic connections are negative.
We derive conditions for fragmentary and essential asymptotic stability of
the three cycles and of the network.

 We prove that one of the cycles in the network  is always completely unstable. One of the other two cycles is essentially asymptotically stable whenever it is fragmentarily asymptotically stable. The third cycle may be fragmentarily asymptotically stable without being essentially asymptotically stable. We also show that at most one of the cycles is fragmentarily asymptotically stable.
This is a necessary condition to guarantee that the whole network is fragmentarily asymptotically stable.
Finally,  we derive conditions for the essential asymptotic stability of the network.
That it is not asymptotically stable follows from our result concerning stability of generic compact robust heteroclinic networks.

\section{Background}

Consider $\Gamma$-equivariant vector fields in $\R^n$. If the vector field is represented by an ordinary differential equation $\dot{x}=f(x)$ then
for all element, $\gamma$, of the compact Lie group $\Gamma$ and for every element, $x$, in $\R^n$ we have
$$
f(\gamma .x)=\gamma .f(x).
$$
The vector field possesses a {\em heteroclinic cycle} if there exist equilibria $\xi_j$, $j=1, \hdots, m$, and trajectories $\kappa_{j-1,j}=[\xi_{j-1} \rightarrow \xi_j]$ for the vector field such that
$$
\kappa_{j-1,j}\subset W^u(\xi_{j-1}) \cap W^s(\xi_j) \neq \varnothing,
$$
where  $\xi_{m+1}=\gamma \xi_1$ for some $\gamma\in\Gamma$.
 In an equivariant context, we identify equilibria and connections in the same group orbit. That is, equilibria $\xi_i$ and $\xi_j$ such that $\xi_i=\gamma \xi_j$ for some $\gamma \in \Gamma$  and connections $\kappa_{j-1,j}=[\xi_{j-1} \rightarrow \xi_j]$ and $\gamma \kappa_{j-1,j}=[\gamma\xi_{j-1} \rightarrow \gamma\xi_j]$ are thought of as the same. 
A {\em heteroclinic network} is a connected set that is the union of two or more heteroclinic cycles.
 Note that in an equivariant context, the Guckenheimer-Holmes example \cite{GuckenheimerHolmes} is a cycle, not a network.

Even though in general heteroclinic connections in cycles are not robust, in the   symmetric context some invariant spaces arise naturally. If restricted to these spaces the connections are from saddle to sink, this ensures robustness of heteroclinic cycles and networks. A {\em fixed-point space} for a subgroup $\Sigma$ of $\Gamma$ is defined as
$$
\mbox{Fix}(\Sigma) = \{ x \in \R^n: \; x = \sigma .x, \;\; \mbox{for all } \;\; \sigma \in \Sigma \}.
$$
We denote by $L_j=$Fix$(\Delta_j)$ the  fixed-point space containing $\xi_j$ and by $P_{ij}=$Fix$(\Sigma_{ij})$ the  fixed-point space containing the heteroclinic connection $\kappa_{ij}$.
In this paper we assume that $L_j$ is 1-dimensional and  that $P_{ij}$ is 2-dimensional.

The dynamics near heteroclinic cycles and networks depends on the stability of the heteroclinic objects. The study of this stability  relies, as usual, on the properties of return maps which are compositions of local and  global maps.
Local maps near  equilibria $\xi_j$ depend on the eigenvalues of the linearisation $df(\xi_j)$.
Local and global maps also depend on the {\em isotypic decomposition} of the complement to $L_j$ and to $P_{ij}$ in $\R^n$
under the actions of $\Delta_j$ and $\Sigma_{ij}$, respectively.
The isotypic decomposition of a space is the unique decomposition into a direct sum of subspaces each of which is the sum of all equivalent irreducible representations.
Here it is  used both to provide the geometric structure  of the global maps
as well as to describe the eigenvalues and eigenspaces at equilibria
(see \cite{Gol&al88} for more detail).

The stability  properties of a heteroclinic cycle or network range from {\em asymptotic stability} (a.s.), the strongest, to {\em complete instability} (c.u.), the weakest, and are defined below.

For a compact invariant set $X\subset\R^n$, and a flow $\Phi_t(x)$,
the $\delta$-basin of attraction of $X$ is
$$
{\cal B}_\delta(X) = \{x\in \R^n~:~d(\Phi_t(x),X)<\delta\hbox{~for all~}t>0\mbox{~and~} \lim_{t\rightarrow+\infty}d(\Phi_t(x),X)=0 \}.
$$
Analogously, the $\delta$-basin of attraction of $X$ for a map is obtained by replacing $t$ by $n$ and $\Phi_t(x)$ by $f^n(x)$ in the set above. The following definitions of stability are relevant in our work.
The concepts in Definitions~\ref{cu} and \ref{def:eas} are from Melbourne \cite{Melbourne1991} while Definition~\ref{fas} is from Podvigina \cite{op12}.
In what follows, $\ell(\cdot)$ denotes the Lebesgue measure in the appropriate context and dimension.

\begin{figure}
\centerline{
\includegraphics[scale=1]{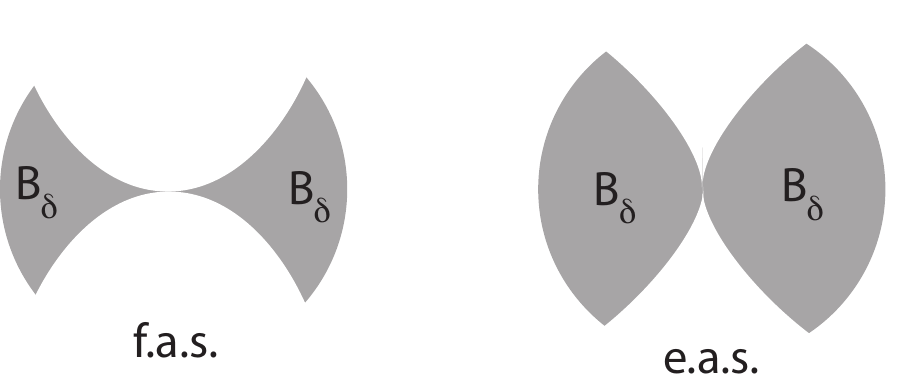}}
\caption{
 Shape of the basins of attraction $ {\cal B}_\delta(X)$ for
the stability types of  Definitions~\ref{fas} (left)  and \ref{def:eas} (right).
The grey region represents
 the intersection of a basin with a cross-section to $X$.}
\label{figStab}
\end{figure}

\begin{definition}\label{cu}
A compact invariant set $X$ is {\it completely unstable (c.u.)} if there exists $\delta>0$
such that $\ell({\cal B}_\delta(X))=0$.
\end{definition}
\begin{definition} \label{fas}
A set $X$ is {\it fragmentarily asymptotically stable (f.a.s.)} if
$\ell({\cal B}_\delta(X))>0$ for any $\delta>0$ (see  Figure~\ref{figStab}).
\end{definition}
\begin{definition} \label{def:eas}
A set $X$ is {\it essentially asymptotically stable (e.a.s)} if
$$
\lim_{\delta\to0}\left(
\lim_{\varepsilon\to0}
\frac{\ell(N_{\varepsilon}(X)\setminus{\cal B}_\delta(X))}{\ell \left(N_{\varepsilon}(X)\right)}
\right)=0,
$$
where $N_{\varepsilon}(X)$ denotes $\varepsilon$-neighbourhood of $X$ (see  Figure~\ref{figStab}).
\end{definition}
\begin{definition}
A set $X$ is {\it  asymptotically stable (a.s)} if for any $\delta>0$ there exists $\varepsilon>0$ such that
$N_{\varepsilon}(X)\subset{\cal B}_\delta(X)$.
\end{definition}
Note that a.s. implies e.a.s., and e.a.s. implies f.a.s. but the converse does not hold:
if a set is a.s. it attracts a full neighbourhood of points;
if a set is e.a.s. it attracts a subset of (asymptotically) full measure in its
neighbourhood; if a set is f.a.s., it attracts a positive measure set from its any neighbourhood.
From the point of view of
simulations and applications, sets that are either a.s. or e.a.s. are the ones
more likely to be observed.

For a compact invariant set $X\subset\R^n$, a point $x\in X$, $\delta>0$ and $ N_\varepsilon(x)$  the ball of centre $x$ and radius   $\varepsilon>0$, let
$$
S_{\varepsilon,\delta}(x)=\frac{ \ell\left({\cal B}_\delta(X)\cap N_\varepsilon(x) \right)}{\ell\left( N_\varepsilon(x) \right)}.
$$

\begin{definition}[\cite{pa11}]\label{ai}
Let $X$ be  a compact invariant set with  $x\in X$.
Define:
$$
\sigma_{loc,-}(x)=\lim_{\delta\to 0}\lim_{\varepsilon\to 0}\frac{\ln S_{\varepsilon,\delta}(x)}{\ln \varepsilon}
\quad\mbox{and}\quad
\sigma_{loc,+}(x)=\lim_{\delta\to 0}\lim_{\varepsilon\to 0}\frac{\ln \left(1-S_{\varepsilon,\delta}(x) \right)}{\ln \varepsilon}
$$
with the conventions that $\sigma_{loc,-}(x)=\infty$ if there is an $\varepsilon_0$ such that $S_{\varepsilon,\delta}(x)=0$ for all
$\varepsilon<\varepsilon_0$, and that
$\sigma_{loc,+}(x)=\infty$ if there is an $\varepsilon_0$ such that $S_{\varepsilon,\delta}(x)=1$ for all
$\varepsilon<\varepsilon_0$.

The  (local) {\em stability index} of $X$ at $x$
is then
$$
\sigma(x,X)=\sigma_{loc,+}(x)-\sigma_{loc,-}(x).
$$
Note that $\sigma_{loc,\pm}\ge 0$, hence $\sigma(x,X)\in[-\infty,\infty]$.
\end{definition}

The stability index $\sigma(x,X)$ is constant for $x$ in a trajectory \cite[Theorem 2.2]{pa11}.
If  $X$ is either a heteroclinic cycle or a compact heteroclinic network having a connection $\kappa_{ij}$, then this allows us to define $\sigma(\kappa_{ij},X)$ as $\sigma(x,X)$, for some  $x\in\kappa_{ij}$.

\section{Stability results}\label{stability_results}

We divide our stability results into two types: those that study the network as a whole and those that study stability of fixed points of maps. The main result concerning stability of a network is of a negative kind. We show that many heteroclinic networks never are asymptotically stable. The results pertaining to fixed points of maps may be applicable to other cycles or networks beyond the present case study.

\subsection{Stability of networks}

In this short subsection, we prove generic results that apply to  robust heteroclinic networks in an Euclidean space of any finite dimension. We provide sufficient conditions that prevent a heteroclinic network in $\R^n$ from being a.s. In particular, we immediately conclude that the network in the case study of this article is not a.s.

\begin{theorem}\label{lem5}
Let $X\subset\R^n$ be a robust heteroclinic cycle or network   with equilibria $\xi_j$.
 Assume that $X$ is compact.
If there exists $\xi_j\in X$ such that $W^u(\xi_j)\not\subset X$ then  $X$ is
not asymptotically stable.
\end{theorem}

\proof
Since $X$ is compact and $W^u(\xi_j)\not\subset X$, there exist $y\in
W^u(\xi_j)$ and $\delta>0$ such that $d(y,X)>\delta$. Denote by $\Phi_t(y)$ the
trajectory through $y$. Since $\lim_{t\to-\infty}\Phi_t(y)=\xi_j$, for any
$\varepsilon>0$ there exists $T_\varepsilon>0$ such that
$y_\varepsilon=\Phi_{-T_\varepsilon}(y)$ satisfies $d(\xi_j,y_\varepsilon)<\varepsilon$.
Hence, for any $\varepsilon>0$ we have
$d(y_\varepsilon,X)<\varepsilon$ and $d(\Phi_{T_\varepsilon}(y_\varepsilon),Y)>\delta$.
\qed

\begin{Cy}\label{onetwo}
Let $X\subset\R^n$ be a compact robust heteroclinic network comprised of
 equilibria $\xi_j$ and a finite number of one-dimensional
connections. Suppose that there exists $\xi_j\in X$ such that
$\dim W^u(\xi_j)\ge 2$. Then $X$ is not asymptotically stable.
\end{Cy}

\proof
Since $X$ is comprised of a finite number of one-dimensional connections,
we have $\dim X=1$. Hence, $W^u(\xi_j)\not\subset X$.
\qed

\begin{Cy}
Let $X\subset\R^n$ be a compact robust heteroclinic network  with equilibria $\xi_j$.
If for some equilibrium $\xi_j$ there is a transverse eigenvalue with positive real part
then $X$ is not asymptotically stable.
\end{Cy}

An example of asymptotically stable heteroclinic network with
a continuum of one-dimensional connections can be found in \cite{kirk2010}.
Apparently Corollary~\ref{onetwo} holds true if a countable number of connections
is assumed, but the current statement assuming a finite number of connections
is sufficient for our purposes.
We remark that an extension of the above results to networks whose nodes
are periodic orbits should be possible. However, networks with more complex
nodes need extra care. These are outside the scope of this article.

Concerning weaker notions of stability, it follows from the definition of f.a.s. that if $X$ is a robust heteroclinic network such that at least one of its cycles is f.a.s. then $X$ is f.a.s. Examples in \cite{cl2014} show that the same does not hold for e.a.s.

\subsection{Stability of fixed points}

The following are technical results useful for the study in Section~\ref{CaseStudy}. We provide conditions for different types of stability of fixed points of maps. These maps take several forms which are common in return maps to cross-sections to connections of heteroclinic cycles or networks.

\begin{lemma}\label{lem1}
Consider the map $\bh(p,q)=(p^{\gamma\alpha}q^{\gamma\beta},p^{\alpha}q^{\beta})$,
$\bh:~\R_+^2\to\R_+^2$.
The fixed point $(p,q)=(0,0)$ of the map $\bh$ is
\begin{itemize}
\item [(i)] f.a.s. if and only if
\begin{equation}\label{lem1i}
\gamma\alpha+\beta>1\hbox{ and }\gamma>0;
\end{equation}
\item [(ii)] e.a.s. if and only if (\ref{lem1i})
and $|\max\{\alpha,\beta\}|>|\min\{\alpha,\beta\}|$;
\item [(iii)] a.s. if and only if (\ref{lem1i}), $\alpha>0$ and $\beta>0$.
\end{itemize}
\end{lemma}

\proof
(i) For $n\ge1$ the iterates $(p_n,q_n)=\bh^n(p_0,q_0)$ satisfy
$p_n=q_n^{\gamma}$ and $q_{n+1}=q_n^{\gamma\alpha+\beta}$, therefore
conditions (\ref{lem1i}) are necessary. To show that the conditions
are sufficient, we note that (\ref{lem1i}) implies that at least one of
$\alpha$ and $\beta$ is positive.
Denote $Q_{\delta}=(0,\delta)^2$. The points $(p_0,q_0)\in Q_{\delta}$
such that
\begin{equation}\label{pt1}
p_0^{\gamma_0\alpha}q_0^{\gamma_0\beta}<\delta,
\end{equation}
where $\gamma_0=\min\{1,\gamma\}$, satisfy $(p_n,q_n)\in Q_{\delta}$
for any $n\ge 0$. Since (\ref{pt1}) is equivalent to
$$p_0^{\alpha}q_0^{\beta}<\delta^{1/\gamma_0},$$
and at least one of $\alpha$ and $\beta$ is positive,
the set of such points has positive measure for any $\delta>0$.

(ii, iii) Since e.a.s. and a.s. imply f.a.s., we assume that the conditions
(\ref{lem1i}) are satisfied. As we noted above, at least one of $\alpha$
and $\beta$ is positive. Let $\alpha>0$.
 From (\ref{pt1}), if $\beta$ is positive then all $(p_0,q_0)\in Q_{\varepsilon}$,
where $0<\varepsilon<\delta^{1/\gamma_0(\alpha+\beta)}$, satisfy
$(p_n,q_n)\in Q_{\delta}$ for any $n\ge 0$, which implies that
the origin is a.s. and e.a.s.

For negative $\beta$ we decompose
$Q_{\delta}=Q_{\delta}^I\cup Q_{\delta}^{II}$, where
$$Q_{\delta}^I=\{(p,q)\in Q_{\delta}~:~
p^{\gamma_0\alpha}q^{\gamma_0\beta}<\delta\},\
Q_{\delta}^{II}=\{(p,q)\in Q_{\delta}~:~
p^{\gamma_0\alpha}q^{\gamma_0\beta}>\delta\}.$$
By construction, $\bh^n(p_0,q_0)\in Q_{\delta}$ for any $n\ge 0$ and
$(p_0,q_0)\in Q_{\delta}^I$, while
$\bh(p_0,q_0)\not\in Q_{\delta}$ for any $(p_0,q_0)\in Q_{\delta}^{II}$.
Since for any $\varepsilon>0$ the set
$Q_{\varepsilon}\cap Q_{\delta}^{II}$ is not empty, the origin is not a.s.

If $\alpha>-\beta$, then
$$\lim_{\varepsilon\to0}
\frac{\ell(Q_{\delta}^{II}\cap Q_{\varepsilon})}{\ell \left(Q_{\varepsilon}\right)}=
\lim_{\varepsilon\to 0}
\frac{-\beta\delta^{-1/\gamma_0\beta}}{\alpha-\beta}\varepsilon^{-\alpha/\beta-1}=0,$$
while for $\alpha<-\beta$
$$\lim_{\varepsilon\to 0}
\frac{\ell(Q_{\delta}^{I}\cap Q_{\varepsilon})}{\ell\left( Q_{\varepsilon}\right)}=
\lim_{\varepsilon\to0}
\frac{-\alpha\delta^{-1/\gamma_0\alpha}}{\alpha-\beta}\varepsilon^{-\beta/\alpha-1}=0.$$
Therefore, for positive $\alpha$ the condition for e.a.s. is that
$\alpha>-\beta$. Similarly, for positive $\beta$ the condition for
e.a.s. is that $\beta>-\alpha$. Both conditions are satisfied if
$|\max\{\alpha,\beta\}|>|\min\{\alpha,\beta\}|$.
\qed

\begin{lemma}\label{lem3}
Consider the matrix
\begin{equation}\label{mtA}
A=\left(
\begin{array}{cc}
\alpha_1&\beta_1\\
\alpha_2&\beta_2
\end{array}
\right).
\end{equation}
If
\begin{equation}\label{condl3}
\alpha_1\ge0,\quad \alpha_2>0\quad\hbox{ and }\quad\det A<0
\end{equation}
then the matrix has real eigenvalues, $\lambda_+>0$ and $\lambda_-<0$, and
$v_{11}v_{12}>0$, where $(v_{11},v_{12})$ is the eigenvector associated with
$\lambda_+$. Furthermore, $|\lambda_+|>|\lambda_-|$ if and only if, additionally,
\begin{equation}\label{condl3ex}
\alpha_1+\beta_2>0.
\end{equation}
\end{lemma}

\proof
Since $\det A<0$, the matrix has one positive real eigenvalue and one negative,
that we denote by $\lambda_+$ and $\lambda_-$,
respectively. Decompose $(1,0)=\bv_1+\bv_2$, where
$\bv_j=(v_{j1},v_{j2})$, $j=1,2$,
$A\bv_1=\lambda_+\bv_1$ and $A\bv_2=\lambda_-\bv_2$.
From
$$\left(\begin{array}{c}1\\0\end{array}\right)=\bv_1+\bv_2\hbox{ and }
A\left(\begin{array}{c}1\\0\end{array}\right)=
\lambda_+\bv_1+\lambda_-\bv_2$$
we obtain that $v_{11}=(\alpha_1-\lambda_-)/(\lambda_+-\lambda_-)>0$ and
$v_{12}=\alpha_2/(\lambda_+-\lambda_-)>0$.

Since $\lambda_++\lambda_-=\alpha_1+\beta_2$, the inequality
$|\lambda_+|>|\lambda_-|$ is satisfied if and only if
$\alpha_1+\beta_2$ is positive.
\qed

\begin{Cy}\label{cor1}
Suppose the matrix $A$ in \eqref{mtA} satisfies \eqref{condl3} and \eqref{condl3ex} and let
$(x_1,y_1)=A(x_0,y_0)$,  where $x_1,y_1,x_0,y_0$ are negative,
$y_1/x_1<y_0/x_0$.
Then the set $V=\{(x,y)\in\R^2_-\ :\ y_1/x_1\le y/x \le y_0/x_0\}$
is $A$-invariant, \emph{i.e.} $AV\subset V$.
\end{Cy}

\proof
Let $(a,b)$ be the coordinates of points in $\R^2$ in the basis comprised
of eigenvectors of the matrix $A$, $\bv^+$ and $\bv^-$.
Denote by $(a_0,b_0)$ and $(a_1,b_1)$ the coordinates of $(x_0,y_0)$ and
$(x_1,y_1)$, respectively, and choose the directions of the eigenvectors
such that $a_0>0$ and $b_0>0$. In the coordinates $(a,b)$ the set $V$ is
$$V=\{(a,b)\in\R^2~~:~~a>0,~~b_1/a_1\le b/a\le b_0/a_0~~\}.$$
Since $\lambda_+^2/\lambda_-^2>1$ and $A(a,b)=(\lambda_+a,\lambda_-b)$,
for any $(a,b)\in V$ we have
$$b_1/a_1=\lambda_-b_0/\lambda_+a_0\le \lambda_-b/\lambda_+a\le
\lambda_-b_1/\lambda_+a_1=\lambda_-^2b_0/\lambda_+^2a_0<b_0/a_0,$$
which implies that $V$ is $A$-invariant.
\qed

\begin{lemma}\label{lem2}
Consider the map $\bh:~\R_+^2\to\R_+^2$,
\begin{equation}\label{maphh}
\bh(p,q)=(p^{\alpha_1}q^{\beta_1},p^{\alpha_2}q^{\beta_2})\hbox{, where }\alpha_2>0.
\end{equation}
The fixed point $(p,q)=(0,0)$ of the map $\bh$ is
\begin{itemize}
\item [(i)] f.a.s. if and only if all the following conditions hold:
\begin{enumerate}
\item $\alpha_1+\beta_2>0$,
\item either $\alpha_1+\beta_2+\beta_1\alpha_2-\alpha_1\beta_2>1$
or $\alpha_1+\beta_2>2$,
\item either $\beta_1>0$ or $\alpha_1-\beta_2>0$;
\item $\left(\alpha_1-\beta_2\right)^2+4\beta_1\alpha_2\ge 0$;
\end{enumerate}
\item [(ii)] a.s. if and only if both conditions below hold:
\begin{enumerate}
\item $\alpha_1>0, \beta_1>0, \beta_2>0$,
\item either $\alpha_1+\beta_2+\beta_1\alpha_2-\alpha_1\beta_2>1$
or $\alpha_1+\beta_2>2$.
\end{enumerate}
\end{itemize}
Moreover,  in case ({\sl i}), if either  $\beta_1>0$ or $\alpha_1+\beta_2+\beta_1\alpha_2-\alpha_1\beta_2>1$ then condition {\sl 4.} is redundant.
\end{lemma}

\proof
Consider the transition matrix $A$ given  in \eqref{mtA} and let
$\lambda_{\max}$ be its  eigenvalue maximal in absolute value, with
associated eigenvector ${\bf v}^{\max}=(v_1^{\max},v_2^{\max})$.
As proved in \cite{op12}, the fixed point is f.a.s. if and only if
$$\lambda_{\max}\hbox{ is real, }\lambda_{\max}>1\hbox{ and }
v_1^{\max}v_2^{\max}>0.$$

The eigenvalues of $A$ are
\begin{equation}\label{eig}
\lambda_{\pm}=\frac{\alpha_1+\beta_2}{2}\pm
\left(\frac{(\alpha_1-\beta_2)^2}{4}+\beta_1\alpha_2\right)^{1/2},
\end{equation}
with the associated eigenvectors
\begin{equation}\label{eigv}
{\bf v}^{\pm}=\left(\frac{\lambda_{\pm}-\beta_2}{\alpha_2},1\right).
\end{equation}
From (\ref{eig}), the eigenvalues are real if and only if
$(\alpha_1-\beta_2)^2+4\beta_1\alpha_2\ge 0$. We have
$|\lambda_+|>|\lambda_-|$ if and only if $\alpha_1+\beta_2>0$. The inequality $\lambda_+>1$ is
satisfied if and only if either $\alpha_1+\beta_2+\beta_1\alpha_2-\alpha_1\beta_2>1$ or $\alpha_1+\beta_2>2$.  Finally,
(\ref{eigv}) implies that
$v_1^{\max}v_2^{\max}=v_1^+v_2^+=(\lambda_+-\beta_2)/\alpha_2>0$ if and only if
either $\beta_1\alpha_2>0$ or $\alpha_1-\beta_2>0$.

Because $\alpha_2>0$ it follows that
$\beta_1>0$ is equivalent to $\beta_1\alpha_2>0$ and this implies that
$(\alpha_1-\beta_2)^2+4\beta_1\alpha_2>0$  when the first part of {\sl 3.} holds.
When $\alpha_1+\beta_2+\beta_1\alpha_2-\alpha_1\beta_2>1$ we can write
\begin{eqnarray*}
\left(\alpha_1-\beta_2\right)^2+4\beta_1\alpha_2&=&
\left(\alpha_1+\beta_2\right)^2-4\left(\alpha_1\beta_2-\beta_1\alpha_2\right)>\\
&>&\left(\alpha_1+\beta_2\right)^2+4\left(1-(\alpha_1+\beta_2)\right)
=\left(\alpha_1+\beta_2-2\right)^2>0
\end{eqnarray*}

\medskip\noindent
Recall \cite{op12} that the fixed point $(0,0)$ of the map (\ref{maphh}) is
a.s. if and only if all entries of the matrix $A$ are non-negative and
$\lambda_{\max}>1$.
Condition {\sl 1.} of part ({\sl ii}) is equivalent to the non-negativity of the entries of $A$.
 Due to (\ref{eig}), $\alpha_1+\beta_2>0$ implies
$\lambda_{\max}=\lambda_+$, and from the arguments above, part (ii) is proven.
\qed

The next two lemmas provide conditions for the stability of the fixed point of a map of the form
$\bh(p,q)=(\max\{p^{\gamma\alpha_2}q^{\gamma\beta_2},
p^{\alpha_1}q^{\beta_1}\},p^{\alpha_2}q^{\beta_2})$, depending on
relations among parameters.
These lemmas are used to study the stability of
cycles ${\cal C}_{123}$ and ${\cal C}_{143}$
in Section~\ref{CaseStudy}, for which the signs are as given
in the statement of the lemmas. So, we prove the lemmas in the restricted form
that is sufficient for our purposes,
although similar proofs may be given for other parameter ranges.

\begin{lemma}\label{lem4new}
Consider the map $\bh:~\R_+^2\to\R_+^2$,
$$\bh(p,q)=(\max\{p^{\gamma\alpha_2}q^{\gamma\beta_2},
p^{\alpha_1}q^{\beta_1}\},p^{\alpha_2}q^{\beta_2})\hbox{, where }$$
$$\alpha_1\ge0,\ \alpha_2>0,\ \gamma\alpha_2-\alpha_1>0\hbox{ and }
\gamma_1=\frac{\beta_1-\gamma\beta_2}{\gamma\alpha_2-\alpha_1}>\gamma.$$
The fixed point $(p,q)=(0,0)$ of the map $\bh$ is
\begin{itemize}
\item not f.a.s. if either $\gamma<0$ or  $\gamma\alpha_2+\beta_2<1;$

\item not e.a.s. if $ \alpha_2<-\beta_2;$

\item not a.s. if $ \beta_2<0$.

\item  f.a.s. if  $\gamma>0$ and  $\gamma\alpha_2+\beta_2>1;$

\item e.a.s. if $\gamma>0$, $\gamma\alpha_2+\beta_2>1$ and $ \alpha_2>-\beta_2;$

\item  a.s. if $\gamma>0$, $\gamma\alpha_2+\beta_2>1$ and $ \beta_2>0.$

\end{itemize}
\end{lemma}

\begin{figure}
\centerline{
\includegraphics[scale=1]{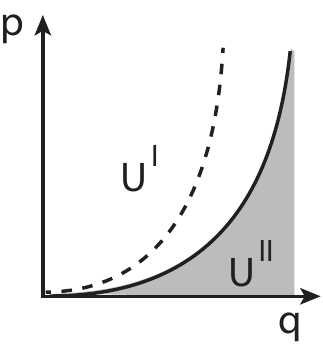}}
\caption{
Under the conditions of Lemma~\ref{lem4new}  with $\gamma>0$
the white  set $U^{I}$ with boundary  in the solid line $p=q^{\gamma_1} $ is mapped by $\bh$ into the dashed curve $p=q^\gamma$, that is contained in $U^I$.
The grey set $U^{II}$ is mapped inside $U^{I}$.Here we show the case $\gamma>1$.}
\label{setUlem4new}
\end{figure}

\proof
Evidently, $\gamma\le 0$ implies that the map is completely unstable, hence till the end of the proof we assume that $\gamma>0$.
The result is local, therefore we work on
$V_{\varepsilon}=\{(p,q)\in \R^2~:~ 0<p<\varepsilon, ~ 0<q<\varepsilon\}$,
with $0<\varepsilon<1$,
that we decompose as $V_{\varepsilon}=U^{I}\cup U^{II}$ where
\begin{equation}\label{eqU11}
U^{I}=\{(p,q)\in V_{\varepsilon}~:~~\bh(p,q)=\bh^{I}(p,q)\equiv
(p^{\gamma\alpha_2}q^{\gamma\beta_2},p^{\alpha_2}q^{\beta_2})\},
\end{equation}
\begin{equation}\label{eqU21}
U^{II}=\{(p,q)\in V_{\varepsilon}~:~~\bh(p,q)=\bh^{II}(p,q)\equiv(p^{\alpha_1}q^{\beta_1},p^{\alpha_2}q^{\beta_2})\}.
\end{equation}
These sets, shown in Figure~\ref{setUlem4new} can also be written as
$$
U^{I}=\{(p,q)\in V_{\varepsilon}~:~~ p\ge q^{\gamma_1}\}
\qquad\mbox{and}\qquad
U^{II}=\{(p,q)\in V_{\varepsilon}~:~~ p\le q^{\gamma_1}\} .
$$
Recall, from the proof of Lemma~\ref{lem1}, that $\bh^{I}(p,q)=(p_1,q_1)=(q^\gamma_1, q_1)$.
Hence, for any $(p,q)\in U^{I}$,  we have
$\bh(p,q)=(p_1,q_1)$ where $p_1= q^\gamma_1> q^{\gamma_1}_1$.
Therefore,  $\bh(U^{I})$ is contained in  the curve $p=q^\gamma$ and in particular
 $\bh(U^{I})\subset U^{I}$ (see Figure~\ref{setUlem4new}).

For  $(p,q)\in U^{II}$,   again let
$\bh(p,q)=(p_1,q_1)=\bh^{II}(p,q)$.
We have that $q_1=p^{\alpha_2}q^{\beta_2}$ and, by definition of $U^{II}$,  that
$p_1=p^{\alpha_1}q^{\beta_1}>p^{\gamma \alpha_2}q^{\gamma\beta_2}=q_1^\gamma>q_1^{\gamma_1}$.
Hence, $\bh(p,q)=(p_1,q_1)\in U^{I}$.

Thus, the conditions for stability are those  given in Lemma~\ref{lem1} for the map $\bh^{I}$.
\qed

\begin{lemma}\label{lem4}
Consider the map $\bh:~\R_+^2\to\R_+^2$,
$$\bh(p,q)=(\max\{p^{\gamma\alpha_2}q^{\gamma\beta_2},
p^{\alpha_1}q^{\beta_1}\},p^{\alpha_2}q^{\beta_2})\hbox{, where }$$
$$
\alpha_1\ge 0,\ \alpha_2>0,\
\ \beta_1\alpha_2-\alpha_1\beta_2>0,
\gamma\alpha_2-\alpha_1>0\hbox{ and }
\gamma_1=\frac{\beta_1-\gamma\beta_2}{\gamma\alpha_2-\alpha_1}<\gamma.
$$

\begin{itemize}
\item[(a)] Assume in addition that $\alpha_1+\beta_2>0$. The fixed point $(p,q)=(0,0)$ of the map $\bh$ is:
\begin{itemize}
\item[\textbullet] not f.a.s. if either $\gamma<0$ or $\alpha_1+\beta_2+\beta_1\alpha_2-\alpha_1\beta_2<1$;
\item[\textbullet] not a.s. if either $\beta_1<0$ or $\beta_2<0$;
\item[\textbullet] f.a.s. if  $\gamma>0$ and $\alpha_1+\beta_2+\beta_1\alpha_2-\alpha_1\beta_2>1$;
\item [\textbullet]a.s. if $\beta_1>0$, $\beta_2>0$,  $\gamma>0$ and $\alpha_1+\beta_2+\beta_1\alpha_2-\alpha_1\beta_2>1$.
\end{itemize}
\item[(b)] Assume in addition that $\alpha_1+\beta_2<0$. The fixed point $(p,q)=(0,0)$ of the map $\bh$ is:
\begin{itemize}
\item[\textbullet]  not f.a.s. if either $\gamma<0$ or $\alpha_2(\gamma\alpha_1+\beta_1)+\beta_2(\gamma\alpha_2+\beta_2)<1$;
\item[\textbullet]  f.a.s. if  $\gamma>0$ and $\alpha_2(\gamma\alpha_1+\beta_1)+\beta_2(\gamma\alpha_2+\beta_2)>1$;
\item[\textbullet] never a.s.
\end{itemize}
\end{itemize}
\end{lemma}

The conditions may be interpreted in terms of the matrix $A$ of \eqref{mtA} and its eigenvalues, as follows:
 $\alpha_1+\beta_2$ is the trace of the matrix $A$, hence it has  the same sign as $|\lambda_+|-|\lambda_-|$.
 Thus $\alpha_1+\beta_2>0$ means that the eigenvalues of $A$ satisfy $|\lambda_+|>|\lambda_-|$.
If $p(\lambda)$ is the characteristic polynomial of $A$, then $p(1)<0$ implies that $\lambda_+>1$ and
$p(1)=-(\alpha_1+\beta_2+\beta_1\alpha_2-\alpha_1\beta_2)+1$.
The condition $\beta_1\alpha_2-\alpha_1\beta_2>0$ means $\det A<0$.

\begin{figure}[hhh]
\centerline{
\includegraphics[scale=1]{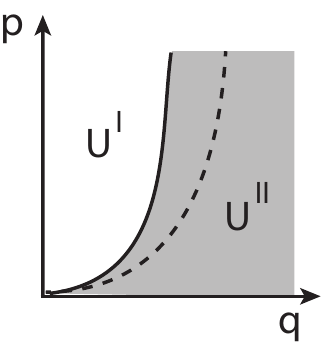}
\qquad\qquad\qquad
\includegraphics[scale=1]{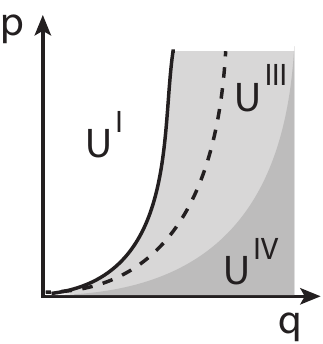}
\qquad\qquad\qquad
\includegraphics[scale=1]{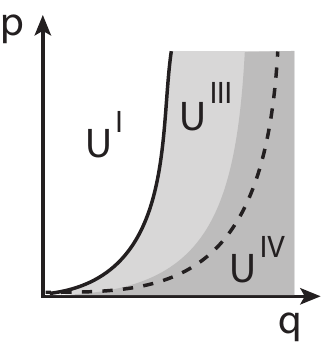}
}
\caption{Under the conditions of Lemma~\ref{lem4}  with $\gamma>0$
the white  set $U^{I}$ with boundary in the solid line
$p=q^{\gamma_1} $ is mapped by $\bh$ into the dashed curve $p=q^\gamma$, that is   contained in  the grey set $U^{II}$.
The latter set can be further subdivided in two components that, in case {\sl (a)} when $\alpha_1+\beta_2>0$, are mapped as
$\bh(U^{IV})\subset U^I$,  $\bh(U^{I})\subset U^{III}$ and $\bh(U^{III})\subset U^{III}$, shown in the middle. Further iterates of points in $U^{III}$ approach a line in $U^{III}$.
 Case{\sl (b)} is shown on the right, where $\bh(U^{IV})\subset U^{I}$, $\bh(U^{I})\subset U^{IV}$ and
for $(p,q)\in U^{III}$ the iterates $\bh^n(p,q)$ either
escape from $V_{\varepsilon}$ or $\bh^n(p,q)\in U^{IV}$ for some
$0<n<\infty$.
\label{setUlem4}}
\end{figure}

\proof
As in the proof of Lemma~\ref{lem4new}, we
decompose $V_{\varepsilon}=U^{I}\cup U^{II}$, where the sets $U^{I}$ and $U^{II}$
are defined in \eqref{eqU11} and \eqref{eqU21}.
 Again, $\bh$ maps $U^I$ into the curve $p=q^\gamma$, but now this curve is contained in $U^{II}$ (see Figure~\ref{setUlem4}).
Let $\gamma_0=(\beta_1-\gamma_1\beta_2)/(\gamma_1\alpha_2-\alpha_1)$,
if $\gamma_0>0$ then
decompose further $U^{II}=U^{III}\cup U^{IV}$, where
$$
U^{III}=\{(p,q)\in U^{II}:~ \bh(p,q)\in U^{II}\}=
\{(p,q)\in U^{II}:~ q^{\gamma_0}\le p\le q^{\gamma_1} \}
$$
and
$$
U^{IV}=\{(p,q)\in U^{II}:~ \bh(p,q)\in U^{I}\}=
\{(p,q)\in U^{II}:~  p\le q^{\gamma_0} \}.
$$
If $\gamma_0\le 0$ then the set $U^{IV}$ is empty and $U^{III}$ coincides with  $U^{II}$.

Let $W_\varepsilon=\{(x,y)\in\R^2_-:~ x<\ln\varepsilon,~ y<\ln\varepsilon \}$.
In the variables $(x,y)=(\ln p,\ln q)$ the sets $U^I$, $U^{III}$ and
$U^{IV}$ are (see Figure~\ref{setUlem4}):
$$
\begin{array}{l}
U^I=\{(x,y)\in W_\varepsilon:~ x\ge\gamma_1y\},\\
U^{III}=\{(x,y)\in W_\varepsilon:~ \gamma_0y\le x\le\gamma_1y\},\\
U^{IV}=\{(x,y)\in W_\varepsilon:~ x\le\gamma_0y\}.
\end{array}
$$
Let $(x_0,y_0)=-(\gamma_0,1)$, $(x_1,y_1)=A(x_0,y_0)$ and
$(x_2,y_2)=A^2(x_0,y_0)$.
Our choice of  $\gamma_0$ and $\gamma_1$ implies
that $x_1/y_1=\gamma_1$ and $x_2/y_2=\gamma$.
In the variables $(a,b)$ employed in the proof of  Corollary~\ref{cor1}
the sets
$U^I$, $U^{III}$ and $U^{IV}$ satisfy:
$$
\begin{array}{l}
U^{III}\subset\{~(a,b)\in\R^2~:~a>0,~b_1/a_1\le b/a\le b_0/a_0~\}\\
U^I\subset\{~(a,b)\in\R^2~:~ b/a\le b_1/a_1~\}\\
U^{IV}\subset\{~(a,b)\in\R^2~:~ b/a\ge b_0/a_0~\}.
\end{array}
$$
As stated, the conditions for stability depend on the
sign of $\alpha_1+\beta_2$. Below we consider the cases of positive
and negative $\alpha_1+\beta_2$ separately. (Note, that generically the sum
does not vanish.)

\medskip
{\sl (a)}
Assume that $\alpha_1+\beta_2>0$. Corollary~\ref{cor1} implies that the
set $U^{III}$ is $\bh$-invariant. In particular,
$(\gamma y',y')=A(\gamma_1y,y)\subset U^{III}$.
For $(a',b')=A(a,b)$, where
$(a,b)\in U^{IV}$ and $a>0$, we have
$b'/a'=\lambda_-b/\lambda_+a<\lambda_-b_0/\lambda_+a_0=b_1/a_1$. Therefore,
$$\left(\bh\left( U^{IV}\right)\cap W_\varepsilon\right)\subset U^{I}.$$

If $(x,y)\in U^I$ then $(x',y')=\bh(x,y)$ satisfies $x'=\gamma y'$.
Therefore, $\left(\bh\left( U^I\right)\cap W_\varepsilon\right)\subset U^{III}$
 (see Figure~\ref{setUlem4}).
In the original variables $(p,q)$ the inclusions can be summarised as
$$
\bh \left(U^{III}\right)\subset U^{III},\quad
 \left(\bh  \left(U^{I}\right)\cap V_\varepsilon\right) \subset U^{III},\quad
 \left(\bh  \left(U^{IV}\right)\cap V_{\varepsilon}\right)\subset U^{I}.
$$
Hence,
$$
\left( \bh\left( \bh\left( V_{\varepsilon} \right)\cap V_{\varepsilon}\right)\cap V_{\varepsilon}  \right) \subset U^{III}
$$
and for any $(p,q)\in V_{\varepsilon}$  and  $n\ge2$  we have
 that either $\bh^n(p,q)\in U^{III}$ or   $\bh^n(p,q)\notin V_{\varepsilon}$.
Therefore, the conditions for
stability are those given in Lemma~\ref{lem2} for the map $\bh^{II}$.

\medskip
{\sl (b)} Assume that $\alpha_1+\beta_2<0$. Therefore, $\beta_2<0$ and the map is not
a.s. To find conditions for f.a.s., note that in the coordinates
$(a,b)$ the iterates $(a_n,b_n)=A^n(a_0,b_0)$ satisfy
$b_n/a_n=\lambda_-^nb_0/\lambda_+^na_0$. Since $|\lambda_-|>|\lambda_+|$,
for any $(a,b)\in U^{III}$ with $b\ne0$ the iterates
$\bh^n(a,b)$ escape from $U^{III}$ for some finite $n>0$.
Moreover, we have $\left(\bh \left(U^{IV}\right)\cap W_\varepsilon\right)\subset U^{I}$ (by the same arguments
as in the case $\alpha_1+\beta_2>0$) and $\left(\bh \left(U^I\right)\cap W_\varepsilon\right)\subset U^{VI}$
(since $(\gamma y',y')=A(\gamma_1y,y)\subset U^{IV}$).
Returning  to the original coordinates $(p,q)$
(see Figure~\ref{setUlem4})
 we proved that
$$
\renewcommand{\arraystretch}{1.5}
\begin{array}{l}
\bh^{n_0}(p,q)\cap V_{\varepsilon}\in U^{IV}
\hbox{for some $n_0>0$, for almost all }(p,q)\in U^{III},\\
\left(\bh \left(U^{IV}\right)\cap V_{\varepsilon}\right)\subset U^I,\qquad
\left(\bh \left(U^I\right)\cap V_{\varepsilon}\right)\subset U^{IV}.
\end{array}
$$
Hence, for almost all initial conditions there exists $n>0$ such that
the iterates $(p_{n},q_{n})=\bh^{n}(p,q)$ satisfy either
$(p_{n},q_{n})\not\in V_{\varepsilon}$ or $(p_{n},q_{n})\in U^{IV}$. In
the latter case $(p_{n},q_{n})=(s^{\gamma},s)$. Further iterates satisfy
$$(p_{n+2},q_{n+2})=
\bh^{n+2}(p_n,q_n)=\bh^2(s^{\gamma},s)=
\bh^I(\bh^{II}(s^{\gamma},s))=$$
$$\bh^I(s^{\alpha_1\gamma+\beta_1},s^{\alpha_2\gamma+\beta_2})=
(s^{\gamma\alpha_2(\alpha_1\gamma+\beta_1)+
\gamma\beta_2(\alpha_2\gamma+\beta_2)},
s^{\alpha_2(\alpha_1\gamma+\beta_1)+\beta_2(\alpha_2\gamma+\beta_2)}),$$
which implies that the map is f.a.s. whenever
$\alpha_2(\alpha_1\gamma+\beta_1)+\beta_2(\alpha_2\gamma+\beta_2)>1$.
\qed

\section{Case study: a heteroclinic network from a convection problem}
\label{CaseStudy}

We study a heteroclinic network supported by the following vector field, given as equations (21) in \cite{clp}:
\begin{equation} \label{equationsR6}
\left\{
\begin{array}{l}
\dot{x}_1 = x_1[\lambda_1+A_1x_1^2+A_2(x_2^2+x_3^2)+C_1y_1^2+C_2(y_2^2+y_3^2)] + A_3x_1x_2^2x_3^2+C_3x_1y_2y_3 \\
\dot{x}_2 = x_2[\lambda_1+A_1x_2^2+A_2(x_1^2+x_3^2)+C_1y_2^2+C_2(y_1^2+y_3^2)] + A_3x_1^2x_2x_3^2+C_3x_2y_1y_3 \\
\dot{x}_3 = x_3[\lambda_1+A_1x_3^2+A_2(x_1^2+x_2^2)+C_1y_3^2+C_2(y_1^2+y_2^2)] + A_3x_1^2x_2^2x_3+C_3x_3y_1y_2 \\
\dot{y}_1 = y_1[\lambda_2+B_1y_1^2+B_2(y_2^2+y_3^2)+C_4x_1^2+C_5(x_2^2+x_3^2)] + B_3y_1y_2^2y_3^2+C_6(y_2x_3^2+y_3x_2^2) \\
\dot{y}_2 = y_2[\lambda_2+B_1y_2^2+B_2(y_1^2+y_3^2)+C_4x_2^2+C_5(x_1^2+x_3^2)] + B_3y_1^2y_2y_3^2+C_6(y_3x_1^2+y_1x_3^2) \\
\dot{y}_3 = y_3[\lambda_2+B_1y_3^2+B_2(y_1^2+y_2^2)+C_4x_3^2+C_5(x_1^2+x_2^2)] + B_3y_1^2y_2^2y_3+C_6(y_1x_2^2+y_2x_1^2) .
\end{array}
\right.
\end{equation}

This vector field is equivariant under the action in $\R^6$ of the   group
$\D_3 \times \Z_2\times  \Z_2$, generated by:
\begin{eqnarray*}
\rho . (x_1,x_2,x_3;y_1,y_2,y_3) & = & (x_2,x_3,x_1;y_2,y_3,y_1) \\
s_1  . (x_1,x_2,x_3;y_1,y_2,y_3) & = & (x_1,x_3,x_2;y_1,y_3,y_2) \\
r  . (x_1,x_2,x_3;y_1,y_2,y_3) & = & -(x_1,x_2,x_3;y_1,y_2,y_3) \\
\gamma^1_{\pi} . (x_1,x_2,x_3;y_1,y_2,y_3) & = &
(x_1,-x_2,-x_3;y_1,y_2,y_3) \\
\gamma^2_{\pi} . (x_1,x_2,x_3;y_1,y_2,y_3) & = &
(-x_1,-x_2,x_3;y_1,y_2,y_3) .
\end{eqnarray*}

\subsection{Description}\label{subDescription}

\begin{figure}
\centerline{
\includegraphics[scale=.8]{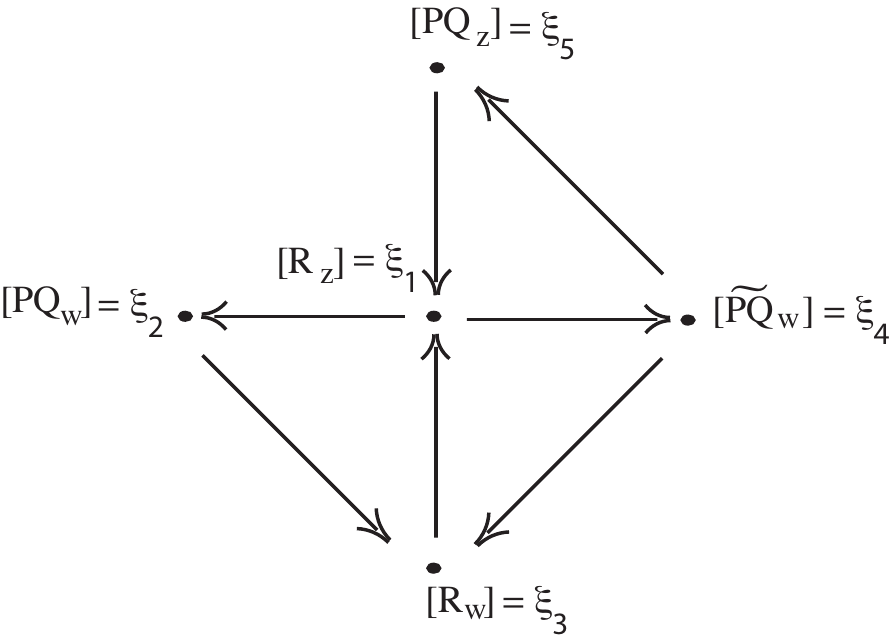}
}
\caption{The network of the case study and its cycles. \label{FigureNetwork}}
\end{figure}

In this subsection we introduce notation that is used in the paper.
The network involves five (isotropy types of) steady states of \eqref{equationsR6} and is shown in Figure~\ref{FigureNetwork}. Here we denote equilibria by $\xi_j$. The correspondence to the notation in  \cite{clp} is as follows:
$$
\begin{array}{lll}
\rz:  & \xi_1=& (x,0,0;0,0,0)\\
\rho^2\qw: & \xi_2=& (0,0,0;0,y,y)\\
\rw:  & \xi_3=& (0,0,0;y,0,0)\\
\rho^2\qtilw& \xi_4=& (0,0,0;0,y,-y)\\
\rho^2\qz: & \xi_5=& (0,x,x;0,0,0)
\end{array}
$$
and $[\xi_j]$ denotes the group orbit of $\xi_j$. By $\kappa_{ij}$ we
denote a heteroclinic connection from $\xi_i$ to $\xi_j$.

By ${\cal C}_{123}$, ${\cal C}_{143}$ and ${\cal C}_{145}$ we denote the
cycles $[\xi_1\to\xi_2\to\xi_3]$, $[\xi_1\to\xi_4\to\xi_3]$ and
$[\xi_1\to\xi_4\to\xi_5]$, respectively.
Conditions for the existence of the network are given in Table 4 of \cite{clp}.

A local basis near $\xi_j$ is comprised of ${\be}_{jk}$, $k=1,\ldots,6$,
which are eigenvectors of $df(\xi_j)$. When $df(\xi_j)$ has an eigenspace
of dimension larger than one, we can use another basis, which is
denoted $\tilde {\be}_{jk}$. For the node $\xi_3$ we will also  need an extra basis, $\hat\be_{3k}$, near $\rho \xi_3$.
The local bases are shown in Table~\ref{tbBases} in Appendix~\ref{AppTable}.
Local  coordinates (i.e. with an origin at $\xi_j$)  at these bases are
denoted by $u_{jk}$, or $\tilde u_{jk}$, respectively.The eigenvalue of
$df(\xi_j)$ associated with ${\be}_{jk}$ is $\lambda_{jk}$.

Tables~\ref{tbIsotypic} and \ref{tbIsotypic2}  of Appendix~\ref{AppTable} provide the  isotypic
 decompositions, while Table~\ref{tab12w} has the eigenvalues and eigenvectors at the nodes.

\subsection{Stability of the cycles}\label{stab}

In this section we derive conditions for stability of individual cycles and the whole network.

\subsubsection{The cycle ${\cal C}_{123}=[\rz\to\qw\to\rw\to\rz]$}\label{cyc1}

In this section we study a cycle that in the quotient space is
$[\rz\to\qw\to\rw\to\rz]$, namely, we derive conditions for the stability of
this cycle. Since the cycle is a part of a network, it is not asymptotically
stable. As shown in \cite{clp}, a trajectory near the cycle
follows equilibria in a certain order.
For definiteness we study asymptotic stability of the cycle
\begin{equation}\label{cyc23}
\rw\to\rz\to\rho^2\qw\to\rho\rw
\qquad\mbox{that is}\qquad
\xi_3\to\xi_1\to\xi_2\to\rho\xi_3 .
\end{equation}
Existence of a trajectory that follows such a sequence of equilibria is shown in \cite[Section 6]{clp} 
and a numerical simulation appears as \cite[Figure 4]{clp}.

Since the cycle is a part of a network and equilibria are stable in the
transverse directions, the eigenvalues $\lambda_{ij}$ of $df(\xi_i)$, $i=1,2,3$, given in Table~\ref{tab12w}, satisfy:
\begin{equation}\label{eig23}
\begin{array}{lcl}
\lambda_{1j}<0,\ j=1,2,3,4,&\qquad& \lambda_{2j}<0,\ j=1,2,3,4,5,\\
\lambda_{3j}<0,\ j=2,4,5,&\qquad&
\lambda_{15},\lambda_{16},\lambda_{26},\lambda_{31}>0.
\end{array}\end{equation}
We remark
that the relative magnitude of $\lambda_{15}$ and $\lambda_{16}$ determines the relative size of the set of points that follow from $\xi_1$ to $\xi_2$ or to $\xi_4$.
If $\lambda_{15}>\lambda_{16}$ then more points follow to $\xi_2$ along this cycle.

We prove the following theorem:
\begin{theorem}\label{th1}
Consider the cycle ${\cal C}_{123}$ and assume that the conditions \eqref{eig23}
are satisfied. Denote
\begin{equation}\label{albt}
\beta_1=\frac{\lambda_{21}\lambda_{35}}{\lambda_{26}\lambda_{31}}(1-\frac{\lambda_{16}}{\lambda_{15}}),\qquad
\beta_2=-\frac{\lambda_{32}}{\lambda_{31}}+
\frac{\lambda_{35}\lambda_{12}}{\lambda_{31}\lambda_{15}}+
\frac{\lambda_{35}\lambda_{22}}{\lambda_{31}\lambda_{26}} \left(1-\frac{\lambda_{16}}{\lambda_{15}}\right)
\end{equation}
Then
\begin{itemize}
\item [(i)] If
$$\lambda_{15}<\lambda_{16}\hbox{ or }\beta_1+\beta_2<1,$$
then the cycle c.u.
\item [(ii)]
If
\begin{equation}\label{condst}
\lambda_{15}>\lambda_{16}\hbox{ and }\beta_1+\beta_2>1,
\end{equation}
then the cycle is e.a.s. The stability indices are:
$$\sigma(\kappa_{31},{\cal C}_{123})=1-\lambda_{16}/\lambda_{15},\quad
\sigma(\kappa_{23},{\cal C}_{123})=+\infty,\quad
\sigma(\kappa_{12},{\cal C}_{123})=+\infty,
$$
\end{itemize}
\end{theorem}

\proof
We approximate the behaviour of trajectories near
the cycle by a return map, which is a composition of local
(approximating behaviour of trajectories near the steady state) and global (approximating
behaviour near heteroclinic connections) maps.
We start by calculating expressions
for  these maps, from which we derive  the expression for the return map
$\tilde \bg:H_3^{\ini}\to H_3^{\ini}$ where $H_3^{\ini}$
is a cross-section of the connection $\rho^2\qw\to\rho\rw$ near $\rho\rw$.
Then we derive conditions for asymptotic stability of the map $\tilde \bg$.
Because of its complexity, the proof of stability conditions of the map $\tilde \bg$
is given in Appendix~\ref{appB}.
Finally, we prove that the cycle is e.a.s. whenever the map $\tilde \bg$ is a.s.
and calculate stability indices of the cycle.

The return map is the superposition of
local maps $\phi_j:\ H^\ini_j\to H^\ou_j$ and global maps
$\psi_{ij}:\ H^\ou_i\to H^\ini_j$.
Here  $H_j^\ini$ and $H_j^\ou $ denote
the cross-sections near $\xi_j$ to the connections to and from $\xi_j$, respectively\footnote{In Section~\ref{netw}, where we deal with the network as a whole,  we  use a more cumbersome notation for the cross-sections, so as to specify the connection.
We also use there the notation $\phi_{312}$ for $\phi_1$ above, to emphasise the connections that are being followed. Since there is no ambiguity here, we use the simpler notation.}.
Cross-sections are taken to be 4-dimensional, since we can disregard the radial direction at each equilibrium.
When we need to specify that the norm of points in the cross-section is smaller than $\varepsilon$ we write $H^\ini_j(\varepsilon)$.
In each cross-section near $\xi_j$, we use coordinates $u_{ji}^\ini$ and $u_{ji}^\ou$ in the direction
of the connections from and to $\xi_i$, respectively.

A local map near $\xi_j\in L_j$,
where $L_j=\Fix\,\Delta_j$, depends on the symmetry group $\Delta_j$, or to be
more precise on the isotypic decomposition of $\R^6\ominus L_j$ under
$\Delta_j$, and on eigenvalues of $df(\xi_j)$.
The isotypic decomposition of $\R^6\ominus L_j$ is
given in Table~\ref{tbIsotypic} of Appendix~\ref{AppTable} and local bases near equilibria
are given in Table~\ref{tbBases} of Appendix~\ref{AppTable}.

The local maps $H^\ini_j\to H^\ou_j$ are obtained from the flow of the linearised equations, as follows.
We compute the flight time from $H^\ini_j$ to $H^\ou_j$ and then $\phi_j$ is obtained
substituting this flight time in the other coordinates, to get:
\begin{equation}\label{phi}
\renewcommand{\arraystretch}{1.5}
\begin{array}{cl}
\phi_1:\ H^\ini_1\to H^\ou_1 &
u_{12}^\ou=u_{12}^\ini|u_{15}^\ini|^{-\lambda_{12}/\lambda_{15}},\
u_{13}^\ou=u_{13}^\ini|u_{15}^\ini|^{-\lambda_{13}/\lambda_{15}},\\
&u_{14}^\ou=D_2|u_{15}^\ini|^{-\lambda_{14}/\lambda_{15}},\
u_{16}^\ou=u_{16}^\ini|u_{15}^\ini|^{-\lambda_{16}/\lambda_{15}}\\
\phi_2:\ H^\ini_2\to H^\ou_2 &
u_{21}^\ou=D_3|u_{26}^\ini|^{-\lambda_{21}/\lambda_{26}},\
u_{22}^\ou=u_{22}^\ini|u_{26}^\ini|^{-\lambda_{22}/\lambda_{26}},\\
&u_{23}^\ou=u_{23}^\ini|u_{26}^\ini|^{-\lambda_{23}/\lambda_{26}},\
u_{24}^\ou=u_{24}^\ini|u_{26}^\ini|^{-\lambda_{24}/\lambda_{26}},\\
\phi_3:\ H^\ini_3\to H^\ou_3 &
u_{32}^\ou=u_{32}^\ini|u_{31}^\ini|^{-\lambda_{32}/\lambda_{31}},\
u_{33}^\ou=u_{33}^\ini|u_{31}^\ini|^{-\lambda_{33}/\lambda_{31}},\\
&u_{35}^\ou=u_{35}^\ini|u_{31}^\ini|^{-\lambda_{35}/\lambda_{31}},\
u_{36}^\ou=D_1|u_{31}^\ini|^{-\lambda_{36}/\lambda_{31}},\\
\end{array}
\end{equation}
where $D_1$, $D_2$ and $D_3$ are positive.

When an equilibrium $\xi_i$ belongs to several different cycles, the local map near it depends on the cycle chosen, since the transverse directions are different.
However we use  the same notation $\phi_i$ for the different local maps at $\xi_i$ and this should not confuse the reader, since the calculations for each cycle are totally independent and occur in different sections.

A global map along
$\kappa_{ij}=[\xi_i\to\xi_j]$, $\kappa_{ij}\subset P_{ij}$, where
$P_{ij}=\Fix\,\Sigma_{ij}$, is predominantly linear.
In order to study stability, it is essential to determine which coefficients of the linear map vanish.
This, in turn, depends on the isotypic
decomposition of $\R^6\ominus P_{ij}$ under $\Sigma_{ij}$
provided in Appendix~\ref{AppTable}.
For the global maps $\psi_{ij}:H^\ou_i\to H^\ini_j$ we take linear approximations:
\begin{equation}\label{psi}
\renewcommand{\arraystretch}{1.3}
\begin{array}{cl}
\psi_{12}:\ H^\ou_1\to H^\ini_2 & u_{22}^\ini=B_1u_{12}^\ou,\
u_{23}^\ini=B_2u_{13}^\ou,\ u_{24}^\ini=B_3u_{14}^\ou,\ u_{26}^\ini=B_4u_{16}^\ou\\
\psi_{23}:\ H^\ou_2\to H^\ini_3 & \hat u_{31}^\ini=C_1u_{22}^\ou+C_2u_{23}^\ou,\\
&\hat u_{32}^\ini=C_3u_{21}^\ou,\ \hat u_{33}^\ini=C_4u_{22}^\ou+C_5u_{23}^\ou,\
\hat u_{35}^\ini=C_6u_{24}^\ou\\
\psi_{31}:\ H^\ou_3\to H^\ini_1 & u_{12}^\ini=A_1\tilde u_{32}^\ou,\
u_{13}^\ini=A_2\tilde u_{33}^\ou,\
u_{15}^\ini=A_3\tilde u_{35}^\ou,\ u_{16}^\ini=A_4\tilde u_{36}^\ou,
\end{array}
\end{equation}
where $A_j$, $B_j$, $C_3$ and $C_6$ are positive.

Therefore, the return map $\tilde{\bg}:\,H^\ini_3\to H^\ini_3$ is
given by the composition
$$
\tilde{\bg}=\psi_{23}\phi_2\psi_{12}\phi_1\psi_{31}\phi_3,
$$
which also
involves the change of coordinates
\begin{equation}\label{ccord}
\renewcommand{\arraystretch}{1.3}
\begin{array}{cc}
\tilde u_{32}^\ini=(u_{32}^\ini+u_{33}^\ini)/\sqrt{2},&
\tilde u_{33}^\ini=(u_{32}^\ini-u_{33}^\ini)/\sqrt{2},\\
\tilde u_{35}^\ini=(u_{35}^\ini+u_{36}^\ini)/\sqrt{2},&
\tilde u_{36}^\ini=(u_{35}^\ini-u_{36}^\ini)/\sqrt{2}.
\end{array}
\end{equation}

Stability properties of the map $\tilde \bg$
are studied in Appendix~\ref{appB}. It is shown that for almost
all (except for a set of zero measure) points in a neighbourhood $V_\varepsilon\subset\R^4$ for
sufficiently small $\varepsilon$ and large $n$ the asymptotic behaviour
of $\tilde \bg^n$ can be approximated by a map $\bh:\R^2_+\to\R^2_+$,
$\bh(p,q)=(\max\{pq^{\beta_2},
q^{\beta_1}\},pq^{\beta_2})$, where
$p=\max\{|x_1|,|x_2|\}$, $q=\max\{|x_3|,|x_4|\}$
and the $\beta$'s are given by (\ref{albt}).
Unfortunately, we cannot
directly apply results of Lemmas~\ref{lem4new} and \ref{lem4}, since there
are sets where components of $\tilde \bg$ vanish, while the components of $\bh$
are non zero.
Similarly to the proof of Lemmas~\ref{lem4new} and \ref{lem4},
we define sets $U^I-U^{III}$, which are subsets of $V_{\varepsilon}$,
and show that, depending on relations
of parameters of the problem, one of these sets is $\tilde \bg$-invariant
for sufficiently small $\varepsilon$, while almost all points (except
for a set of zero measure) in the
other sets are mapped into the invariant set. In this invariant
set the map $\tilde \bg$ is approximated either by $\bh^I$ or by $\bh^{II}$.

The results  in Appendix~\ref{appB} can be summarised as
follows:
\begin{itemize}
\item [$\bullet$]  Lemma~\ref{lem51} proves that for $|\bx|<\varepsilon$ and sufficiently small $\varepsilon>0$
the map $\tilde \bg(\bx)$ can be approximated as
$$
\begin{array}{l}
\tilde\bg(x_1,x_2,x_3,x_4)\approx\\
\approx(\tilde g_1(x_1,x_2,x_3,0),
\tilde g_2(x_1,x_2,x_3,0),\tilde g_3(x_1,x_2,x_3,0),
-\frac{A_4}{ A_3}\tilde g_3(x_1,x_2,x_3,0)).
\end{array}
$$
\item [$\bullet$] According to  Lemmas~\ref{lem81}--\ref{lem84},  if
$$
\lambda_{15}<\lambda_{16},
\quad\hbox{ or }\quad
\beta_1<0,
\quad\hbox{ or }\quad
\beta_2<0,
\quad\hbox{ or }\quad
\beta_1+\beta_2<1,$$
then the origin is a completely unstable fixed point of the map $\tilde \bg$.
\item [$\bullet$]   In Lemma~\ref{lem7}  we prove that if
$$
\lambda_{15}>\lambda_{16},\quad
\beta_1>0,\quad
 \beta_2>0
 \quad\hbox{ and }\quad
 \beta_1+\beta_2>1,$$
then  the origin is an asymptotically stable fixed point of  the map $\tilde \bg$.
\end{itemize}

Evidently, instability of $\tilde \bg$ implies instability of the cycle,
hence (i) is proven.

In order to prove (ii)  we calculate the stability indices for the heteroclinic connections.
Recall, that the stability index is  constant
along a heteroclinic connection and that it can be calculated on
a codimension one surface transverse to the connection \cite{pa11}.
Moreover, since the equilibria in the cycle are stable in the radial direction, we can further restrict the problem to 4 dimensions.

Under the hypotheses of Lemma~\ref{lem7} (see also \eqref{eig23}, \eqref{condst} and use $\lambda_{15}>\lambda_{16}$) we know that the origin is a.s. for $\tilde{\bg}$.
We start by looking at the connection $\kappa_{23}$: consider $\bx\in H^\ini_3({\varepsilon})$, i.e. $|\bx|<\varepsilon$.
From (\ref{phi}) and  (\ref{psi}) for $\tilde \bg^{(3)}(\bx)\in H^\ini_1$ and
$\tilde \bg^{(1)}\tilde \bg^{(3)}(\bx)\in H^\ini_2$ we obtain that
$$
|\tilde \bg^{(3)}(\bx)|<G_1\varepsilon^{s_1}
\quad\hbox{ and }\quad
|\tilde \bg^{(1)}\tilde \bg^{(3)}(\bx)|<G_2\varepsilon^{s_2},
$$
where $\tilde\bg^{(3)}=\phi_{31}\psi_3$ and $\tilde\bg^{(1)}=\phi_{12}\psi_1$.
Here $s_1$ and $s_2$ depend on $\lambda_{ij}$, $G_j>0$ depend on
constants of the local and global maps and on the eigenvalues.
The inequalities (\ref{eig23}) imply
that $s_1>0$ and $s_2>0$. Moreover, the coordinates of
$\bu=\tilde \bg^{(3)}(\bx)$ satisfy $A_3u_3\approx -A_4u_4$.
Since $\lambda_{15}>\lambda_{16}$, for
any $\delta'>0$ there exists $\varepsilon'>0$ such that for any
$\bx\in H^\ini_3({\varepsilon'})$  the following inequalities hold true:
$$
|\tilde \bg^{(3)}(\bx)|<\delta',\qquad
|\tilde \bg^{(1)}\tilde \bg^{(3)}(\bx)|<\delta'
\quad\hbox{ and }\quad
|\tilde \bg^{(2)}\tilde \bg^{(1)}\tilde \bg^{(3)}(\bx)|<\delta'.
$$
Since the origin is asymptotically stable under the map
$\tilde\bg=\tilde \bg^{(2)}\tilde \bg^{(1)}\tilde \bg^{(3)}$ ,
there exists $\varepsilon>0$ such that
$$\tilde\bg^n(\bx)<\varepsilon'\hbox{ for all }n\ge 0\hbox{ and }|\bx|<\varepsilon.$$
Hence,
$$
\lim_{n\to\infty}\tilde\bg^n(\bx)=0,\quad
\lim_{n\to\infty}\tilde\bg^{(3)}\tilde\bg^n (\bx)=0
\quad\hbox{ and }\quad
\lim_{n\to\infty}\tilde\bg^{(1)}\tilde \bg^{(3)}\tilde\bg^n (\bx)=0$$
and
$$
|\tilde\bg^n(\bx)|<\delta',\quad
 |\tilde\bg^{(3)}\tilde\bg^n(\bx)|<\delta'
\quad \hbox{ and }\quad
|\tilde\bg^{(1)}\tilde \bg^{(3)}\tilde\bg^n(\bx)|<\delta'.
$$
Therefore, at the points in the cross-sections, the distance between a trajectory and the
cycle is bounded by $\delta'$ and vanishes as $n\to\infty$.
Linearity of global maps
implies the existence of a constant $C$ such that, taking $\delta'=C\delta$,
 the distance between the trajectory and the cycle is less than $\delta$.
That is, we proved that $\sigma(\kappa_{23},{\cal C}_{123})=+\infty$.
The proof that $\sigma(\kappa_{12},{\cal C}_{123})=+\infty$ is similar and
we omit it.

For the connection $\kappa_{31}$, note that
in $H^\ini_1$ the trajectories that escape $\delta$-neighbourhood of $\xi_1$
along the connection $\kappa_{12}$ satisfy
$u_{16}u_{15}^{-\lambda_{16}/\lambda_{15}}<\delta$. By the same arguments as
above, all such trajectories stay close to the cycle and are attracted to it
as $t\to\infty$. Hence,
$\sigma(\kappa_{31},{\cal C}_{123})=1-\lambda_{16}/\lambda_{15}$.
When $\lambda_{15}>\lambda_{16}$ all stability indices are positive and by
\cite[Theorem 3.1]{Lohse15} the cycle is e.a.s.
\qed

\subsubsection{The cycle ${\cal C}_{143}=[\rz\to\qtilw\to\rw\to\rz]$}\label{cyc2}

In this section we derive conditions for f.a.s. and
calculate stability indices for a cycle that in the quotient space is
$[\rz\to\qtilw\to\rw\to\rz]$. 
Three numerical simulations of this cycle appear in Figures 5--7 of \cite{clp}. 
We consider behaviour of trajectories near the
cycle that in $\R^6$ is
\begin{equation}\label{cyc143}
\rw\to\rz\to\rho^2\qtilw\to-\rho\rw
\qquad\mbox{that is}\qquad
\xi_3\to\xi_1\to\xi_4\to-\rho\xi_3 .
\end{equation}

Since the cycle is a part of a network and by assumption equilibria are stable in the
transverse directions, the eigenvalues $\lambda_{ij}$  of $df(\xi_i)$, $i=1,3,4$, satisfy:
\begin{equation}\label{eig43}
\begin{array}{ll}
\lambda_{1j}<0,\ j=1,2,3,4,
&
\lambda_{4j}<0,\ j=1,4,6,\\
\lambda_{3j}<0,\ j=2,4,5,
&
\lambda_{15},\lambda_{16},\lambda_{42}=\lambda_{43},\lambda_{45},\lambda_{31}>0.
\end{array}
\end{equation}
We prove the following theorem:
\begin{theorem}\label{th2}
Consider the cycle ${\cal C}_{143}$ and assume that the conditions \eqref{eig43}
are satisfied. Denote
\begin{equation}\label{albt4}
\beta_1=\frac{\lambda_{41}\lambda_{35}}{\lambda_{45}\lambda_{31}}(1-\frac{\lambda_{15}}{\lambda_{16}}),\qquad
\beta_2=-\frac{\lambda_{32}}{\lambda_{31}}+
\frac{\lambda_{35}\lambda_{12}}{\lambda_{31}\lambda_{16}}+
\frac{\lambda_{35}\lambda_{42}}{\lambda_{31}\lambda_{45}} \left(1-\frac{\lambda_{15}}{\lambda_{16}}\right)
\end{equation}
Then
\begin{itemize}
\item [(i)] If
$$\lambda_{16}<\lambda_{15}
\quad\hbox{ or }\quad
\beta_1+\beta_2<1
\quad\hbox{ or }\quad
\beta_2<0,$$
then the cycle c.u.
\item [(ii)]
If
\begin{equation}\label{condf4}
\lambda_{16}>\lambda_{15},\quad
 \beta_1+\beta_2>1
 \quad\hbox{ and }\quad
 \beta_2>0
\end{equation}
then the cycle is f.a.s. The stability indices are:
$$
\sigma(\kappa_{43},{\cal C}_{143})=+\infty,$$

$$
\sigma(\kappa_{14},{\cal C}_{143})=
\bigg\{
\begin{array}{lll}
1-\lambda_{42}/\lambda_{45}&\hbox{ if }& \lambda_{45}>\lambda_{42}\\
\lambda_{45}/\lambda_{42}-1&\hbox{ if }& \lambda_{45}<\lambda_{42}
\end{array}
$$

$$
\sigma(\kappa_{31},{\cal C}_{123})=
\bigg\{
\begin{array}{lll}
\min\{1-\lambda_{15}/\lambda_{16},(1+\beta_3-\beta_4)/\beta_4\}&\hbox{ if }&
\beta_4< \beta_3+1\\
\min\{1-\lambda_{15}/\lambda_{16},(-1-\beta_3+\beta_4)/\beta_3,-\beta_3+\beta_4-1\}&
\hbox{ if }&\beta_4\ge \beta_3+1
\end{array}
$$
where

$$
\beta_3=-\frac{\lambda_{12}}{\lambda_{16}}+
\frac{\lambda_{15}\lambda_{42}}{\lambda_{16}\lambda_{45}},\qquad
\beta_4=-\frac{\lambda_{42}}{\lambda_{45}}.
$$
\end{itemize}
\end{theorem}

We begin the proof of the theorem by proving a lemma.
Denote by $w(h,\varepsilon,a,b)$, where $0<\varepsilon<h$, $a>0$ and $b>0$,
the volume of the set $W\subset\R^3$,
$$W=\{\bx\in\R^3~:~x_1x_2^a>hx_3^b,~0\le x_j\le \varepsilon,\hbox{ for }j=1,2,3~\}.$$
\begin{lemma}\label{lemth2}
For sufficiently small $\varepsilon>0$ the volume of the set $W$ is
$$w(h,\varepsilon,a,b)=
\bigg\{
\begin{array}{lll}
C_1(h,a,b)\varepsilon^{2+(a+1)/b}&\hbox{ if }&b< 1+a\\
C_2(h,a,b)\varepsilon^{2+(b-1)/a}+C_3(h,a,b)\varepsilon^{2+b-a}&\hbox{ if }&b\ge 1+a
\end{array}
$$
where $C_i(h,a,b)$ are positive constants independent of $\varepsilon$.
\end{lemma}

\proof
If $b< 1+a$ then for sufficiently small $\varepsilon$ the function
$g(x_1,x_2)=h^{-1/b}x_1^{1/b}x_2^{a/b}$ satisfies $g(x_1,x_2)<\varepsilon$ for
all $0<x_1<\varepsilon$ and $0<x_2<\varepsilon$. Therefore,
$$
w(h,\varepsilon,a,b)=
\int_0^{\varepsilon}\int_0^{\varepsilon}g(x_1,x_2)\rd x_1\rd x_2=
\frac{h^{-1/b}b^2}{ (a+b)(b+1)}\varepsilon^{2+(1+a)/b}.
$$

For $b\ge 1+a$ we represent $W=W^1\setminus W^2$, where
$$
\begin{array}{l}
W^1=\{\bx\in\R^3~:~0<x_3<g(x_1,x_2),~0\le x_1\le \varepsilon,\
0\le x_2\le \varepsilon~\},\\
W^2=\{\bx\in\R^3~:~\varepsilon<x_3<g(x_1,x_2),~0\le x_1\le \varepsilon,\
0\le x_2\le \varepsilon~\}.
\end{array}
$$
Therefore,
$$w(h,\varepsilon,a,b)=
\int_0^{\varepsilon}\int_0^{\varepsilon}g(x_1,x_2)\rd x_1\rd x_2-
\int_{h^{-1/b}\varepsilon^{(b-1)/a}}^{\varepsilon}
\int_{h\varepsilon^bx_2^{-a}}^{\varepsilon}g(x_1,x_2)\rd x_1\rd x_2=$$
\begin{equation}\label{sum2}
\frac{h^{1/a-1/b}b^2}{ (a+b)(b+1)}\varepsilon^{2+(b-1)/a}+
\frac{hb}{ (1-a)(b+1)}\varepsilon^{2+b-a}-
\frac{h^{1/a}b}{ (1-a)(b+1)}\varepsilon^{2+(b-1)/a}.
\end{equation}
\qed

\begin{remark}\label{rlemth2}
For $h\to0$ in the sum (\ref{sum2}) the third term is asymptotically smaller
than the first one.
\end{remark}

\begin{remark}\label{r2lemth2}
In the limit $\varepsilon\to0$ for $a>1$ in the sum (\ref{sum2}) the first term
is asymptotically larger than the second one, while for $a<1$ the second term
is asymptotically larger.
\end{remark}

\noindent
{\it Proof of the theorem.}
Similarly to the proof of Theorem~\ref{th1}, we approximate the behaviour of
trajectories near the cycle by the return map
$\tilde \bg:H_3^{\ini}\to H_3^{\ini}$. For the cycle
${\cal C}_{143}$ the expression for this map
that we derive coincides (up to expressions for coefficients $\beta_1$ and
$\beta_2$) with the one obtained in Theorem~\ref{th1} for the cycle
${\cal C}_{123}$. Hence, we apply results of Appendix~\ref{appB}
to find conditions for asymptotic stability. Calculation of stability indices
for the cycle ${\cal C}_{143}$ is more difficult, because the equilibrium
$\xi_4$ has a three-dimensional unstable manifold, while the unstable manifold
of $\xi_2$ is one-dimensional.

The local maps $H^\ini_j\to H^\ou_j$ are:
\begin{equation}\label{phi2}
\renewcommand{\arraystretch}{1.5}
\begin{array}{cl}
\phi_1:\ H^\ini_1\to H^\ou_1 &
u_{12}^\ou=u_{12}^\ini|u_{16}^\ini|^{-\lambda_{12}/\lambda_{16}},\
u_{13}^\ou=u_{13}^\ini|u_{16}^\ini|^{-\lambda_{12}/\lambda_{16}},\\
&u_{14}^\ou=D_2|u_{16}^\ini|^{-\lambda_{14}/\lambda_{16}},\
u_{15}^\ou=u_{15}^\ini|u_{16}^\ini|^{-\lambda_{15}/\lambda_{16}}\\
\phi_4:\ H^\ini_4\to H^\ou_4 &
u_{41}^\ou=D_3|u_{45}^\ini|^{-\lambda_{41}/\lambda_{45}},\
u_{42}^\ou=u_{42}^\ini|u_{45}^\ini|^{-\lambda_{42}/\lambda_{45}},\\
&u_{43}^\ou=u_{43}^\ini|u_{45}^\ini|^{-\lambda_{42}/\lambda_{45}},\
u_{44}^\ou=u_{44}^\ini|u_{45}^\ini|^{-\lambda_{44}/\lambda_{45}},\\
\phi_3:\ H^\ini_3\to H^\ou_3 &
u_{32}^\ou=u_{32}^\ini|u_{31}^\ini|^{-\lambda_{32}/\lambda_{31}},\
u_{33}^\ou=u_{33}^\ini|u_{31}^\ini|^{-\lambda_{32}/\lambda_{31}},\\
&u_{35}^\ou=u_{35}^\ini|u_{31}^\ini|^{-\lambda_{35}/\lambda_{31}},\
u_{36}^\ou=D_1|u_{31}^\ini|^{-\lambda_{35}/\lambda_{31}},\\
\end{array}
\end{equation}
where $D_1$, $D_2$ and $D_3$ are some  positive constants.

The global maps $\psi_{ij}:H^\ou_i\to H^\ini_j$ are:
\begin{equation}\label{psi2}
\renewcommand{\arraystretch}{1.5}
\begin{array}{cl}
\psi_{14}:\ H^\ou_1\to H^\ini_4 &
 u_{42}^\ini=B_1u_{12}^\ou+B_2u_{13}^\ou,\
u_{43}^\ini=B_3u_{12}^\ou+B_4u_{13}^\ou,\\
&\ u_{44}^\ini=B_5u_{14}^\ou,\ u_{45}^\ini=B_6u_{15}^\ou\\
\psi_{43}:\ H^\ou_4\to H^\ini_3 &
\hat u_{31}^\ini=C_1u_{42}^\ou+C_2u_{43}^\ou,\
\hat u_{32}^\ini=C_3u_{41}^\ou,\\
& \hat u_{33}^\ini=C_4u_{42}^\ou+C_5u_{43}^\ou,\
\hat u_{35}^\ini=C_6u_{44}^\ou\\
\psi_{31}:\ H^\ou_3\to H^\ini_1 &
 u_{12}^\ini=A_1\tilde u_{32}^\ou,\
u_{13}^\ini=A_2\tilde u_{33}^\ou,\\&
u_{15}^\ini=A_3\tilde u_{35}^\ou,\ u_{16}^\ini=A_4\tilde u_{36}^\ou,
\end{array}
\end{equation}
where $A_j$, $B_j$, $C_3$ and $C_6$ are positive.
To complete the  return map $\tilde\bg:H^\ini_3\to H^\ini_3$ one should apply
the change of coordinates (\ref{ccord}) between $\psi_{43}$ and $\phi_{3}$.

Note the similarity of expressions (\ref{phi2}) and (\ref{psi2}) with
the ones (\ref{phi}) and (\ref{psi}).
Here $\psi_{14}$  differs slightly from
$\psi_{12}$, but this does not modify the final expression for superposition.
Hence, we can apply results of Appendix~\ref{appB} about stability of the map
$\tilde\bg$. For the ${\cal C}_{143}$ cycle the conditions for
stability take the form
\begin{itemize}
\item [$\bullet$] If
$$
\lambda_{16}<\lambda_{15},
\quad\hbox{ or }\quad
\beta_1<0,
\quad\hbox{ or }\quad
\beta_2<0,
\quad\hbox{ or }\quad
\beta_1+\beta_2<1,$$
then the origin is a completely unstable fixed point of the map $\tilde \bg$.
\item [$\bullet$] If
$$
\lambda_{16}>\lambda_{15},\quad
\beta_1>0,\quad
 \beta_2>0
 \quad\hbox{ and }\quad
 \beta_1+\beta_2>1,$$
then  the origin is an asymptotically stable fixed point of  the map $\tilde \bg$.
\end{itemize}
Statement {\em (i)} holds true, because instability of $\tilde \bg$ implies instability of the
cycle. Below we prove  {\em (ii)}. The stability properties of the cycle are
studied by calculating stability indices along the connections,
as it was done in the proof of Theorem~\ref{th1}.

Since $\beta_1>0$, the inequalities (\ref{condf4}) imply that
the origin is an a.s. fixed point of the map $\tilde\bg$.
Consider $\bx\in H^\ini_3({\varepsilon})$, i.e. $|\bx|<\varepsilon$.
For almost all $\bx$ (i.e., except the points that belong to the stable manifolds of
the equilibria), the trajectory $\Phi_t(\bx)$ starting at $\bx$ follows
the connection $\kappa_{31}$ and then $\kappa_{14}$, the latter happens
since $\lambda_{16}>\lambda_{15}$. Then, the trajectory follows the
connection $\kappa_{42}$, because the map $(\psi_{43})^{-1}\tilde\bg:H_3^{\ini}\to H_4^{\ou}$
is a superposition of a linear map and an asymptotically stable $\tilde\bg$.
By the same arguments as employed in the proof of Theorem~\ref{th1},
the trajectory stays close to the cycle for all positive $t$,
hence
$\sigma(\kappa_{43},{\cal C}_{143})=+\infty$.

For $\bu\in H^\ini_4$, $\bu=\left(u_{42},u_{43},u_{44},u_{45}\right)$,
the trajectories $\Phi_t(\bu)$ that escape the $\delta$-neigh\-bour\-hood of $\xi_4$
along the connection $\kappa_{43}$ satisfy
\begin{equation}\label{setst2}
u_{42}u_{45}^{-\lambda_{42}/\lambda_{45}}<\delta\hbox{ and }
u_{43}u_{45}^{-\lambda_{42}/\lambda_{45}}<\delta.
\end{equation}
Then, they are mapped by $\psi_{43}$ to $H_3^{\ini}$ and, given that $|\bu|$
is sufficiently small, stay close to the cycle for all $t>0$, as proven above.
The stability index can be positive or negative, depending on the sign
of $\lambda_{42}-\lambda_{45}$. Calculating the measure of the area bounded
by (\ref{setst2}), we obtain that the index is
$1-\lambda_{42}/\lambda_{45}$ for $\lambda_{45}>\lambda_{42}$ and
$\lambda_{45}/\lambda_{42}-1$ for $\lambda_{45}<\lambda_{42}$.

In $H^\ini_1$ the trajectories that escape a $\delta$-neighbourhood of $\xi_1$
along the connection $\kappa_{14}$ satisfy
$u_{15}u_{16}^{-\lambda_{15}/\lambda_{16}}<\delta$.
By substituting the expressions for $\phi_1$ and $\psi_{14}$
into $\phi_4$ (see (\ref{phi2}) and (\ref{psi2})),
we obtain that trajectories that escape  a $\delta$-neighbourhood of $\xi_4$
along the connection $\kappa_{43}$ satisfy
\begin{equation}\label{setst}
pu_{16}^{\beta_3}u_{15}^{\beta_4}<\delta,
\end{equation}
where $p=\max\{|x_{12}|,|x_{13}|\}$ and $\beta_3$ and
$\beta_4$ are the ones given in the statement of the theorem.
The measure of the set (\ref{setst}) is calculated in Lemma~\ref{lemth2}.
Applying the definition of stability indices, we complete
the proof of  {\em (ii)}.
\qed

\begin{Cy}
Consider the cycle ${\cal C}_{143}$ of Theorem~\ref{th2}.
Then:
\begin{itemize}
\item [(iii)]
If  (\ref{condf4})  holds and in addition the inequality $\lambda_{45}>\lambda_{42}$
holds true then the cycle is e.a.s.
\item [(iv)]
If $\lambda_{45}<\lambda_{42}$ then the cycle is not e.a.s.
\end{itemize}
\end{Cy}
\noindent
{\it Proof.}
To prove {\em (iii)} and {\em (iv)}, we note that the only stability index that can be
non-positive is $\sigma(\kappa_{14},{\cal C}_{143})$.
If $\lambda_{45}>\lambda_{42}$ then the index is positive and, hence, the cycle is e.a.s.
If $\lambda_{45}<\lambda_{42}$ then the index is negative and, hence, the cycle is not e.a.s.
\qed

\subsubsection{The cycle ${\cal C}_{145}=[\rz\to\qtilw\to\qz\to\rz]$}\label{cyc3}

In this section we are concerned with the cycle $[\rz\to\rho^2\qtilw\to\rho^2\qz\to\rho^2\rz]$,
that is $\xi_1\to\xi_4\to\xi_5\to\rho^2\xi_1$.
This cycle is  pseudo-simple (see   \cite[Definition~5]{pc16})
because there is a two-dimensional isotypic component corresponding to  the expanding
eigenspace of $\xi_4$.
For the same reason this cycle is completely unstable, as we show in the next theorem.
In \cite{clp} no numerical simulations show trajectories that stayed close to this cycle.

\begin{theorem}\label{th3}
Generically, the cycle ${\cal C}_{145}$ is completely unstable.
\end{theorem}

\proof
The proof is similar to the proof of Theorem~1 in \cite{pc16}.
We consider the map $\phi_4\psi_{14}\phi_1:H^\ini_1\to H^\ou_4$ where $\phi_4$ and $\phi_1$ are the local maps around $\xi_4$ and $\xi_1$, respectively and $\psi_{14}$ is the global map along the connection $\kappa_{14}$.

The local maps are:
\begin{equation}\label{phi3}
\renewcommand{\arraystretch}{1.5}
\begin{array}{cl}
\phi_1:\ H^\ini_1\to H^\ou_1 &
u_{12}^\ou=D_1|u_{16}^\ini|^{-\lambda_{12}/\lambda_{16}},\
u_{13}^\ou=u_{13}^\ini|u_{16}^\ini|^{-\lambda_{12}/\lambda_{16}},\\
&u_{14}^\ou=u_{14}^\ini|u_{16}^\ini|^{-\lambda_{14}/\lambda_{16}},\
u_{15}^\ou=u_{15}^\ini|u_{16}^\ini|^{-\lambda_{15}/\lambda_{16}}\\
\phi_4:\ H^\ini_4\to H^\ou_4 &
u_{41}^\ou=D_2|u_{42}^\ini|^{-\lambda_{41}/\lambda_{42}},\
u_{43}^\ou=u_{43}^\ini|u_{42}^\ini|^{-\lambda_{43}/\lambda_{42}},\\
&u_{44}^\ou=u_{44}^\ini|u_{42}^\ini|^{-\lambda_{44}/\lambda_{42}},\
u_{45}^\ou=u_{45}^\ini|u_{42}^\ini|^{-\lambda_{45}/\lambda_{42}}
\end{array}
\end{equation}
where $D_1$, $D_2$  are some  positive constants.
The expression for the global map $\psi_{14}$ is given in \eqref{psi2}.

For small values of $\bu$,   the coordinate $u_{13}^\ou$  in the expression for $\phi_1$ is much smaller than $u_{12}^\ou$.
Thus, when computing $\psi_{14}\phi_1$  the second terms  in the sums $B_1u_{12}^\ou+B_2u_{13}^\ou$ and
$B_3u_{12}^\ou+B_4u_{13}^\ou$ may be ignored.
Because $\lambda_{43}=\lambda_{42}$ the coordinate $u_{43}^\ou$ in $\phi_4$ may be rewritten as
 $u_{43}^\ou=u_{43}^\ini|u_{42}^\ini|^{-1}$.
Therefore in the final superposition, one obtains $u_{43}^\ou\approx B_3/B_1$.
Since generically $B_3\ne 0$ the term $u_{43}^\ou$ cannot be made arbitrarily small.
\qed

\subsection{The network}\label{netw}

In this section $\Sigma$  denotes the network, $H_{ij}^\ini$ denotes a cross-section to the connection $\kappa_{ij}$ close to the node $\xi_j$ and $H_{ij}^\ou$ is a cross-section to the connection $\kappa_{ij}$ close to the node $\xi_i$.
We also use the more cumbersome notation $\phi_{ijk}:H_{ij}^\ini\rightarrow H_{jk}^\ou$ for what was denoted $\phi_j$ above, to emphasise the connections $\kappa_{ij}$ and $\kappa_{jk}$ that are being followed.

A necessary condition for existence of the network and its stability in the
transverse directions is that eigenvalues $\lambda_{ij}$ satisfy
\begin{equation}\label{eigall}
\begin{array}{lcl}
\lambda_{1j}<0,\ j=1,2,3,4,&\qquad& \lambda_{2j}<0,\ j=1,2,3,4,5,\\
\lambda_{3j}<0,\ j=2,4,5,&\qquad&\lambda_{4j}<0,\ j=1,4,6,\\
\lambda_{5j}<0,\ j=1,2,4,5,6&\qquad&
\lambda_{15},\lambda_{16},\lambda_{26},\lambda_{31},
\lambda_{42}=\lambda_{43},\lambda_{45},\lambda_{53}>0.
\end{array}\end{equation}

\begin{figure}[hhb]
\begin{center}
\includegraphics[scale=.8]{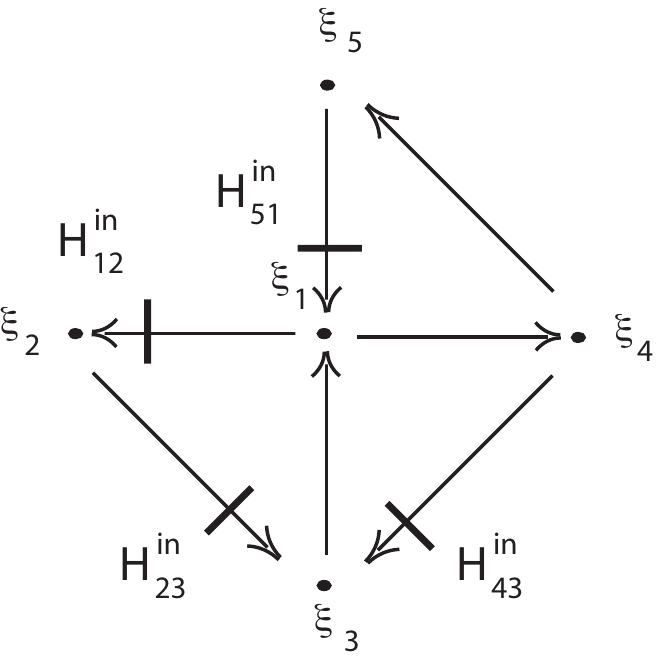}
\qquad\qquad
\includegraphics[scale=.8]{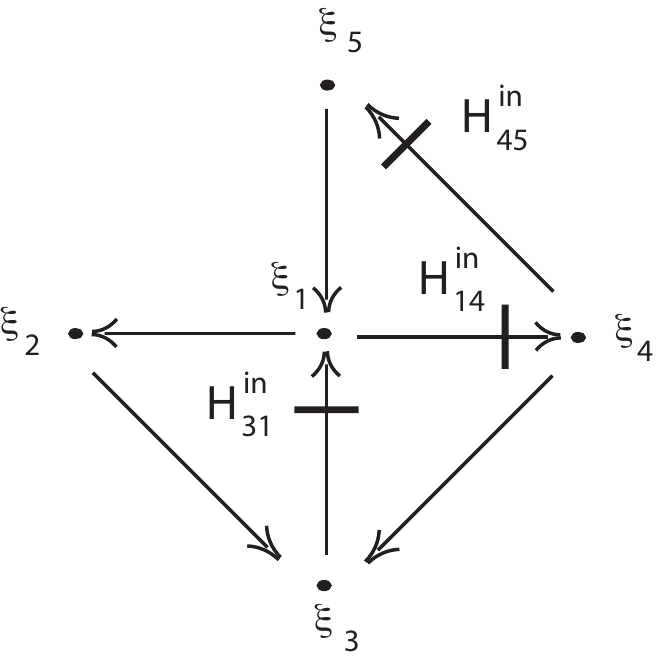}
\end{center}
\caption{
Left: position of the sections in Lemma~\ref{lemmaUnstableSection}. Right: sections in Corollary~\ref{CorolUnstableSection}. }
\label{figSectionsNetwork}
\end{figure}

\begin{lemma}\label{lemmaUnstableSection}
If both cycles ${\cal C}_{123}$ and ${\cal C}_{143}$ are c.u. and $\lambda_{15}>\lambda_{16}$ then
for $\delta>0$ sufficiently small, we have
$\ell \left(H_{ij}^\ini\cap {\cal B}_\delta(\Sigma) \right)=0$,
for $ij\in\{12,23,43,51\}$,
where $\ell$ denotes the Lebesgue measure in $\R^n$ with  $n=\dim H_{ij}^\ini$.
\end{lemma}

\proof
We do the proof for $H_{23}^\ini$, the  cases $H_{12}^\ini$ and $H_{43}^\ini$  are similar.
Let $\bu\in H_{23}^\ini\cap {\cal B}_\delta(\Sigma)$.
Consider the transition map $H_{23}^\ini\to H_{31}^\ini$, denoted by $\tilde{\bg}^{(3)}$ in Lemma~\ref{lem51}.
The estimate \eqref{varrr} shows that $\tilde{\bg}^{(3)}(\bu)$ satisfies   $u_{16}^\ini\approx C u_{15}^\ini$ for some nonzero constant $C$.
If the trajectory of $\bu$ follows the connection $\kappa_{14}$ after passing near $\xi_1$, then from the expression \eqref{phi2} for the local map $\phi_1$ and the estimate above, we get
$$
u_{15}^\ou=u_{15}^\ini|u_{16}^\ini|^{-\lambda_{15}/\lambda_{16}}
\approx
C |u_{15}^\ini|^{1-\lambda_{15}/\lambda_{16}} .
$$
Since $1-\lambda_{15}/\lambda_{16}<0$, we have $\displaystyle \lim_{\delta\to 0} u_{15}^\ou=\infty$.
 Hence, if the trajectory of $\bu$ follows  $\kappa_{14}$ then  $\bu\not\in {\cal B}_\delta(\Sigma)$ for sufficiently small $\delta$.
Therefore, the set $ H_{23}^\ini\cap {\cal B}_\delta(\Sigma)$ is contained in the    set
$ H_{23}^\ini\cap {\cal B}_\delta({\cal C}_{123})$  that, for sufficiently small $\delta$, has zero measure.

For $ H_{51}^\ini$, we  use
$\ell \left(H_{43}^\ini\cap {\cal B}_\delta(\Sigma) \right)=0$ and apply  similar arguments
recalling that by Theorem~\ref{th3} the cycle ${\cal C}_{145}$ is c.u.
\qed

\begin{Cy}\label{CorolUnstableSection}
Under the conditions of Lemma~\ref{lemmaUnstableSection} and
for $\delta>0$ sufficiently small we have
$\ell\left(H_{ij}^\ini\cap {\cal B}_\delta(\Sigma) \right)=0$,
for $ij\in\{45, 14, 31\}$.
\end{Cy}
\proof
Except for a measure zero set, trajectories starting in $H_{45}^\ini$ go to $ H_{51}^\ini$, where we can apply the result of Lemma~\ref{lemmaUnstableSection}.
Similarly, most  trajectories starting in $H_{14}^\ini$ either follow $\kappa_{45}$ or $\kappa_{43}$.
Those following $\kappa_{43}$ end up mostly in $H_{43}^\ini$, and those near $\kappa_{45}$ end up mostly in $H_{45}^\ini$, and both these sets meet ${\cal B}_\delta(\Sigma)$ in a measure zero set.
The arguments for $H_{31}^\ini$ are entirely similar.
\qed

\begin{theorem}\label{th4}
Generically for the network $\Sigma$:
\begin{itemize}
\item [(i)] At most one of the cycles ${\cal C}_{123}$ or ${\cal C}_{143}$ is
f.a.s.
\item [(ii)] The network is f.a.s. whenever one of the cycles is f.a.s.
\end{itemize}
\end{theorem}

\proof

\begin{itemize}
\item [{\em (i)}]
If $\lambda_{16}>\lambda_{15}$ then
Theorem~\ref{th1} implies that if ${\cal C}_{123}$ is c.u. whereas if  $\lambda_{15}>\lambda_{16}$
then ${\cal C}_{143}$ is c.u. by Theorem~\ref{th2}.

\item [ {\em (ii)}]  Clearly, if one of the cycles is f.a.s. then the network is f.a.s.
It remains to see that when both cycles are c.u. then the network is not f.a.s.
If $\lambda_{15}>\lambda_{16}$  this is a consequence of Lemma~\ref{lemmaUnstableSection} and its corollary.

The proof in the case $\lambda_{16}>\lambda_{15}$ is postponed till after we obtain a few lemmas.
\end{itemize}

\begin{figure}
\begin{center}
\includegraphics[scale=.8]{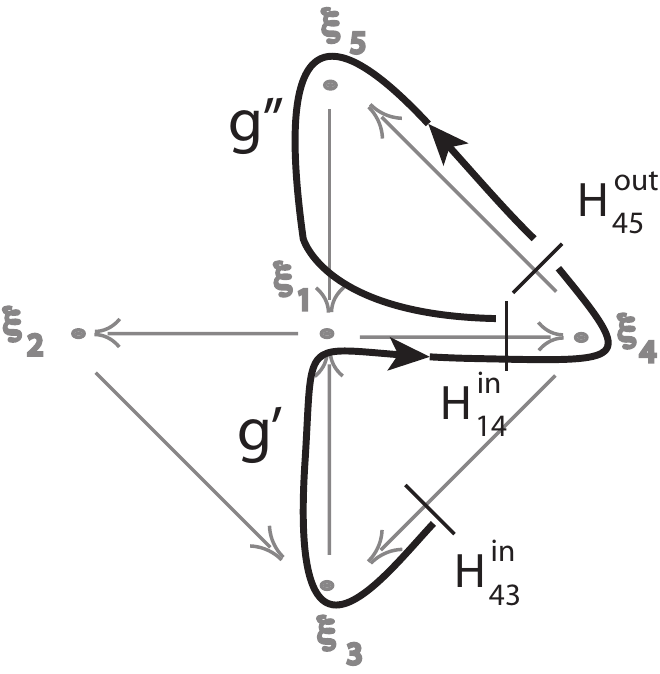}
\end{center}
\caption{
The maps $\bg'$ and $\bg''$ of Lemmas~\ref{nl1} and \ref{nl5}, with the network shown in grey.} 
\label{figLemma17Network}
\end{figure}

\begin{lemma}\label{nl1}
Consider the map $g^\prime:H_{43}^\ini\rightarrow H_{45}^\ou$ given by
$g^\prime=\phi_{145}\psi_{14}\phi_{314}\psi_{31}\phi_{431} $ (see Figure~\ref{figLemma17Network}).
The points in $H_{43}^\ini$ that are mapped by $g^\prime$ into $V_\delta\cap H_{45}^\ou$, for sufficiently small $\delta>0$
belong to the set
$$
V_{45}=\left\{
\bu \in H_{43}^\ini:\ \left(B_3^\prime u_{32}^\ini+B_4^\prime u_{33}^\ini\right)<
\delta c^\prime \max\left(|u_{32}^\ini|,|u_{33}^\ini|\right)
\right\}
$$
where $B_3^\prime$, $B_4^\prime$ and $c^\prime$  are constants, independent on $\delta$.
\end{lemma}

\proof
A direct computation using \eqref{phi3} for $\phi_{145}$ and \eqref{phi2} and \eqref{psi2} for the remaining maps,
shows that writing $g^\prime(\bu)=\left(u_{41}^\ou,u_{42}^\ou,u_{43}^\ou,u_{44}^\ou \right)$
we have
$$
u_{43}^\ou=
 \left(B_1^\prime u_{32}^\ini+B_2^\prime u_{33}^\ini\right)^{-1}
  \left(B_3^\prime u_{32}^\ini+B_4^\prime u_{33}^\ini\right)
$$
\qed

The next lemma follows immediately from  Lemmas~\ref{lem4new} and \ref{lem4}.

\begin{lemma}\label{nl2}
Let $\bh(p,q)$ be one of the maps considered in Lemmas~\ref{lem4new} and \ref{lem4}.
For given  $(p,q)=(r^a,r^b)$  we define $a_n$, $b_n$ by $\bh^n(r^a,r^b)=(r^{a_n},b^{b_n})$.
For any $a>0$, $b>0$ and $s>0$ there exists $0<n_0<\infty$ such that
$\min\{a_{n},b_{n}\}> 0$ for all $n>n_0$ and
at least one of the following is satisfied:
\begin{itemize}
\item[(a)]
$\min\{a_{n_0},b_{n_0}\}\le 0$;
\item[(b)]
$\gamma-s<a_n/b_n<\gamma+s$ for all $n=n^\prime+2m>n_0$;
\item[(c)]
$v_1^+/v_2^+-s<a_n/b_n<v_1^+/v_2^++s$ for all $n>n_0$;
\end{itemize}
where $\gamma$, $v_1^+$ and $v_2^+$ take the meanings they have in  Lemmas~\ref{lem4new} and \ref{lem4}.
\end{lemma}

\begin{lemma}\label{nl3}
Consider the map  $\tilde{\bg}:H_{43}^\ini\rightarrow H_{43}^\ini$  in the proof of Theorem~\ref{th2}.
Let $W\subset H_{43}^\ini$ satisfy
$$
W\subset \left\{ \bu=(u_{31}^\ini=q,u_{32}^\ini, u_{33}^\ini ,u_{35}^\ini):\quad
a_jq^{\gamma_j}<u_{3j}^\ini<b_jq^{\gamma_j}\quad j=2,3
\right\}
$$
where $a_jb_j>0$, $\gamma_j>0$ and $\gamma_2\ne \gamma_3$.
Then, for sufficiently small $\varepsilon>0$,
any point $\hat\bu=\tilde\bg(\bu)=(\hat u_{31}^\ini=\hat q,\hat u_{32}^\ini, \hat u_{33}^\ini ,\hat u_{35}^\ini)\in
\tilde\bg(W\cap V_\varepsilon)$
satisfies
$$
\hat a_j\hat q^{\hat \gamma_j}<\hat u_{3j}^\ini<\hat b_j\hat q^{\hat \gamma_j}\quad j=2,3
$$
where $\hat a_j\hat b_j>0$ are independent on $\varepsilon$, $\hat \gamma_3=1$ and
$\displaystyle \hat\gamma_2=\frac{\beta_1}{\min\{\gamma_2,\gamma_3\} +\beta_2 }$.
\end{lemma}

\proof
Follows from Lemma~\ref{lem51} and from the expressions
\eqref{varrr} and \eqref{mapg0}.
\qed

\begin{Cy}\label{nc1}
Generically, the statement of Lemma~\ref{nl3} holds if $\tilde\bg$ is replaced by $\tilde\bg^n$ for any finite $n>0$, with a different $\hat\gamma_2$.
\end{Cy}

\begin{definition}
We say that a set $V\subset\R^4$ is {\em conical} with exponents $(1,\gamma_2,\gamma_3,\gamma_4)$, $\gamma_j>0$, if
$$
V\subset\left\{ (x_1,x_2,x_3,x_4)\in\R^4:\quad a_j x_1^{\gamma_j}< x_j<b_j x_1^{\gamma_j}, \ j=,2,3,4
\quad\mbox{where}\quad a_jb_j>0
\right\}.
$$
\end{definition}

\begin{lemma}\label{nl4}
Generically all our maps $\bg:H_{ij}^\ini\rightarrow H_{k\ell}^\ini$ have the following property:
if an initial set $U\subset H_{ij}^\ini$ is conical then for sufficiently small $\varepsilon>0$ the image $\bg(U\cap V_\varepsilon)$ is also conical.
\end{lemma}

\proof
The property holds  generically for each one of the local and global maps.
Hence, it also holds for compositions of these maps.
\qed

\begin{lemma}\label{nl5}
Denote $\hat\bg=\bg''\bg' $ where $\bg'':H_{45}^\ou\rightarrow H_{14}^\ini$ is given by
$\bg''=\psi_{14}\phi_{514} \psi_{51}\phi_{451}\psi_{45}$ and $\bg'$ was defined in Lemma~\ref{nl1}
(see Figure~\ref{figLemma17Network}).
Let $U_{45}\subset H_{43}^\ini\cap \cB_\delta(\Sigma) $ be the subset that is mapped to $H_{45}^\ou$ by $\bg'$.
For any $c_0>0$ and sufficiently small $\varepsilon>0$ there is a set $W_{c_0}$ such that
the image
$\hat\bg\left(  (U_{45}\setminus W_{c_0})\cap V_\varepsilon \right)$ is conical with exponents
$(1,1,\lambda_{16}-\lambda_{14}, \lambda_{16}-\lambda_{15})$.
The set $W_{c_0}$ is defined by
$$
W_{c_0}=\left\{
\left| L_{j1} r q^{-1}+L_{j2}p^{\gamma_1}q^{\gamma_2}+L_{j3}p^{\gamma_3}q^{\gamma_4} \right|
<c_0 \max\left\{ |rq^{-1}|, |p^{\gamma_1}q^{\gamma_2}|, |p^{\gamma_3}q^{\gamma_4}|\right\}
\ \  j=1,2,3
\right\}
$$
where $q=u_{31}^\ini$, $p=u_{32}^\ini$, $r=B_3' u_{32}^\ini+B_4' u_{33}^\ini$
and $L_{ji}$,$i, j=1,2,3$, $\gamma_j$, $ j=1,2,3,4$ and  $B_3'$, $B_4' $ are constants coming from the expressions for global and local maps.
\end{lemma}

\proof
Let $\bu\in H_{45}^\ini$ be given by $\bu=(u_{51}^\ini,u_{53}^\ini,u_{54}^\ini,u_{55}^\ini)$.
The coordinates $u_{5j}^\ini$, $j=3,4,5$ for the map $\psi_{45}:H_{45}^\ini:\rightarrow H_{45}^\ou$ have the form
$$
u_{5j}^\ini=
L_{j-2, 1}'u_{43}^\ou+L_{j-2, 2}'u_{44}^\ou+L_{j-2, 3}'u_{45}^\ou
\qquad
j=3,4,5 .
$$
This implies that the composition $\psi_{45} \bg':H_{43}^\ini\rightarrow H_{45}^\ini$, in terms of $p$, $q$, $r$ satisfies
$$
u_{5,j+2}^\ini=
L_{j1} r q^{-1}+L_{j2}p^{\gamma_1}q^{\gamma_2}+L_{j3}p^{\gamma_3}q^{\gamma_4}
\qquad
j=1,2,3 .
$$
If $(p,q,r)$ satisfies
\begin{equation}\label{pqr}
\left| L_{j1} r q^{-1}+L_{j2}p^{\gamma_1}q^{\gamma_2}+L_{j3}p^{\gamma_3}q^{\gamma_4} \right|
\ge c_0 \max\left\{ |rq^{-1}|, |p^{\gamma_1}q^{\gamma_2}|, |p^{\gamma_3}q^{\gamma_4}|\right\}
\quad j=1,2,3
\end{equation}
there exist $a_j$, $b_j$ =, $j=1,2$ such that
$$
a_1|u_{53}^\ini|<|u_{54}^\ini|<b_1|u_{53}^\ini|
\qquad\qquad
a_2|u_{53}^\ini|<|u_{55}^\ini|<b_2|u_{53}^\ini|.
$$
By  writing expressions for $\phi_{415}$ and $\psi_{51}$ and proceeding as above, it can be shown that
if $(p,q,r)$ satisfies \eqref{pqr} then for $\bu\in H_{51}^\ini$  there are $a_j'$, $b_j'$ =, $j=1,2$ such that  we have
$$
a_1'|u_{14}^\ini|<|u_{15}^\ini|<b_1'|u_{14}^\ini|
\qquad\qquad
a_2'|u_{14}^\ini|<|u_{16}^\ini|<b_2'|u_{14}^\ini|.
$$
Using the expressions  \eqref{psi2} for $\phi_{514}$ and \eqref{phi3} for  $\phi_{14}$ we obtain the exponent  \\$(1,1,\lambda_{16}-\lambda_{14}, \lambda_{16}-\lambda_{15})$ of the statement.
\qed
\bigbreak

\noindent{\em End of proof of Theorem~\ref{th4}.}
In the case $\lambda_{16}>\lambda_{15}$,
from Lemmas~\ref{nl1}, \ref{nl2}, \ref{nl3}, \ref{nl4} and Corollary~\ref{nc1} it follows that, for sufficiently small $\varepsilon>0$ almost all trajectories starting in
$$
\hat\bg\left((U_{45}\setminus W_{c_0})\cap V_\varepsilon \right),
$$
after making one turn around ${\cal C}_{145}$,
never again leave the $\delta$-neighbourhood of ${\cal C}_{143}$.
If we take $c_0>0$ small and for small $\varepsilon>0$, we may treat $W_{c_0}$ as in Appendix~\ref{appB}
to show that almost all trajectories that remain close to the network are attracted to the cycle ${\cal C}_{143}$.
Since this cycle is not f.a.s. this implies that the set of these trajectories has zero measure, i.e.
that for small $\delta>0$ we have $\ell\left( H_{43}^\ini\cap\cB_\delta(\Sigma)\right)=0$.

The proof for the other cross-sections follows  arguments similar  to those in Lemma~\ref{lemmaUnstableSection} and its corollary.
\qed

\begin{theorem}\label{th5}
Consider the network $\Sigma$ and assume that the eigenvalues $\lambda_{ij}$
satisfy the inequalities (\ref{eigall}). Then
\begin{itemize}
\item [(i)] suppose the conditions \eqref{condst} hold.
\begin{itemize}
\item[(a)] If $\lambda_{45}<\lambda_{42}$ then $\Sigma$ is not e.a.s.
\item[(b)] If $\lambda_{45}>\lambda_{42}$ then $\Sigma$ is e.a.s.
\end{itemize}
\item [(ii)] suppose  the conditions \eqref{condf4} hold.
\begin{itemize}
\item[(a)] If $\lambda_{45}<\lambda_{42}$ then $\Sigma$ is not e.a.s.
\item[(b)] If
$-\lambda_{12}\lambda_{45}+\lambda_{15}\lambda_{42}<\lambda_{16}\lambda_{42}$
then $\Sigma$ is not e.a.s.
\item[(c)] If $\lambda_{45}>\lambda_{42}$ and
$-\lambda_{12}\lambda_{45}+\lambda_{15}\lambda_{42}>\lambda_{16}\lambda_{42}$
then $\Sigma$ is e.a.s.
\end{itemize}
\end{itemize}
\end{theorem}

\proof
Let $V_\varepsilon$ be an $\varepsilon$-neighbourhood of the origin.
Note that at $\xi_4$ there are three outgoing connections: one to $\xi_3$, one to $\xi_5$, and another to $\gamma_\pi^2\xi_5$.

{\em (ia)} Let $\lambda_{42}>\lambda_{45}$.
Consider the set $Q_{\varepsilon}\subset H_{14}^\ini$ defined as
$Q_{\varepsilon}=\left( H_{14}^\ini\cap\cB_\delta(\Sigma)\right)\cap V_{\varepsilon}$.
The set can be decomposed as the union
$Q_{\varepsilon}=Q_{42}\cup Q_{43}\cup Q_{45}$, such that trajectories from
$Q_{4j}$ leave the $\delta$-neighbourhood of $\xi_4$ along the connection
tangent to $\be_{4j}$. The sets $Q_{4j}$ satisfy
\begin{equation}\label{estset}
\begin{array}{l}
Q_{42}\subset\{~(u_{42}^\ini,u_{43}^\ini,u_{44}^\ini,u_{45}^\ini)~:~
u_{43}^\ini<\delta u_{42}^\ini\hbox{ and }
u_{45}^\ini<\delta( u_{42}^\ini)^{\lambda_{45}/\lambda_{42}}~\}\\
Q_{43}\subset\{~(u_{42}^\ini,u_{43}^\ini,u_{44}^\ini,u_{45}^\ini)~:~
u_{42}^\ini<\delta u_{43}^\ini\hbox{ and }
u_{45}^\ini<\delta( u_{43}^\ini)^{\lambda_{45}/\lambda_{42}}~\}\\
Q_{45}\subset\{~(u_{42}^\ini,u_{43}^\ini,u_{44}^\ini,u_{45}^\ini)~:~
u_{43}^\ini<\delta( u_{45}^\ini)^{\lambda_{42}/\lambda_{45}}\hbox{ and }
u_{42}^\ini<\delta( u_{45}^\ini)^{\lambda_{42}/\lambda_{45}}~\}.
\end{array}
\end{equation}
From inclusions \eqref{estset} we obtain that measures of the sets satisfy
$$\ell(Q_{42})<\frac{1}{2}\delta\varepsilon^4,\quad
\ell(Q_{43})<\frac{1}{2}\delta\varepsilon^4,\quad
\ell(Q_{45})<\frac {\delta\lambda_{42}}{\lambda_{42}+\lambda_{45}}
\varepsilon^{3+\lambda_{42}/\lambda_{45}}.$$
Hence, for sufficiently small $\varepsilon$ we have
$\ell(Q_{\varepsilon})<2\delta\varepsilon^4$, which
according to Definition \ref{def:eas} implies that $\Sigma$ is not e.a.s.

{\em (ib)} Suppose that $\lambda_{42}<\lambda_{45}$.
From  conditions \eqref{condst} it follows
(see the proof of Theorem \ref{th2}) that
$\sigma(\kappa_{12},\Sigma)=\sigma(\kappa_{23},\Sigma)=+\infty$
and $\sigma(\kappa_{12},\Sigma)>0$. Arguments similar to those applied in
the proof of Theorem \ref{th2} imply that $\sigma(\kappa_{43},\Sigma)=+\infty$ and
$\sigma(\kappa_{51},\Sigma)\ge \lambda_{15}/\lambda_{16}-1>0$.
Calculating the measure of the set  $Q_{45}$ constructed above,
we obtain that
$\sigma(\kappa_{14},\Sigma)\ge \lambda_{45}/\lambda_{42}-1>0$.

To estimate $\sigma(\kappa_{45},\Sigma)$, we introduce the set
$W_{\varepsilon}\subset H_{45}^\ini\cap V_{\varepsilon}$ comprised of the points
$\bu=(u_{51}^\ini,u_{53}^\ini,u_{54}^\ini,u_{55}^\ini)\in H_{45}^\ini$
that are mapped by $\phi_{514}\psi_{51}\phi_{451}$ to $H_{12}^\ini$.
The coordinates of these points satisfy
\begin{equation}\label{es45}
\begin{array}{l}
|L'_{31}u_{54}^\ini|u_{53}^\ini|^{-\lambda_{54}/\lambda_{53}}+
L'_{32}u_{55}^\ini|u_{55}^\ini|^{-\lambda_{55}/\lambda_{53}}+
L'_{33}C|u_{55}^\ini|^{-\lambda_{56}/\lambda_{53}}|<\\
\delta|L'_{21}u_{54}^\ini|u_{53}^\ini|^{-\lambda_{54}/\lambda_{53}}+
L'_{22}u_{55}^\ini|u_{55}^\ini|^{-\lambda_{55}/\lambda_{53}}+
L'_{23}C|u_{55}^\ini|^{-\lambda_{56}/\lambda_{53}}|^{\lambda_{16}/\lambda_{15}}
\end{array}
\end{equation}
(Here $L'_{ij}$ and $C$ are the constants of the maps $\psi_{51}$ and
$\phi_{451}$, respectively.) By straightforward but lengthy integration
(similar to the one in the proof of Lemma \ref{lemth2}, but significantly
longer and with bulky final result), it
can be shown that the measure of the set $W_{\varepsilon}$ satisfies
$\ell(W_{\varepsilon})>\varepsilon^4-\delta\varepsilon^{3+s}$, where $s>1$
depends of the exponents $\lambda_{ij}$ that are involved in (\ref{es45}).
Hence,  $\sigma(\kappa_{45},\Sigma)>0$ and part {\em (ib)} is proven.

The proof of {\em (iia)} is identical to the proof of {\em (ia)}. To prove {\em (iib)}, we note
that the points in $H_{51}^\ini$ that are mapped by
$\phi_{341}\psi_{14}\phi_{415}$ to $H_{43}^\ou$ satisfy
$$D|u_{16}^\ini|^{\lambda_{12}/\lambda_{16}}<\delta
|u_{15}^\ini(u_{16}^\ini)^{\lambda_{15}/\lambda_{16}}|^{\lambda_{42}/\lambda_{45}}.$$
As it is shown in the proof of Theorem \ref{th3}, the points in $H_{51}^\ini$
that are mapped
neither to $H_{12}^\ou$ by $\phi_{512}$ nor to $H_{43}^\ou$
by $\phi_{341}\psi_{14}\phi_{415}$, escape from the $\delta$-neghbourhood
of the cycle. Hence,
$$\sigma(\kappa_{52},\Sigma)\le \max\{\lambda_{15}/\lambda_{16}-1,
-\lambda_{12}\lambda_{45}/\lambda_{16}\lambda_{42}+\lambda_{15}/\lambda_{16}-1\}<0.$$
The proof of {\em (iic)} is similar to the proof of  {\em (ib)} and is omitted.
\qed

\section{Conclusion}
We  complete the study of stability of the heteroclinic
network emerging in an ODE obtained from the equations of Boussinesq
convection by the center manifold reduction \cite{clp}.
 We derive and prove
conditions for fragmentary asymptotic stability and essential asymptotic
stability for the network and individual cycles it is comprised of.

This is the first systematic study of stability of heteroclinic cycles that are not of
type A or a generalisation of type Z
(see \cite{GdSC2017} for the latter).
Although we consider a particular  case study, the
proposed approach consisting of well-defined steps is applicable to other
heteroclinic cycles in $\R^n$. Moreover, some of the lemmas that we prove
are not restricted to the case under investigation and can become useful in
other systems.

The study of stability of heteroclinic networks is less common
than that of heteroclinic cycles.
It requires the construction and composition of several transition maps between
 cross-sections to connections belonging to different cycles.
 This procedure is likely to work for other networks as
well.

Our results show that derivation of general stability conditions
for heteroclinic cycles in $\R^n$ with $n\ge6$ is a highly non-trivial task,
if   at all possible. It would be of interest to identify classes
of heteroclinic cycles for which derivation of stability conditions
is possible. To do so, one should somehow classify possible
maps $h$ obtained at step (b) discussed in the introduction.
The classification (at least, partial) should start
with determining possible forms of the map $h$.

In Section 3 we prove Theorem \ref{lem5} stating necessary conditions
for a heteroclinic network to be asymptotically stable.
Corollary \ref{onetwo} of this theorem implies that the network under
investigation in not asymptotically stable.
The theorem can be used to prove instability of more general types of
heteroclinic networks, than considered in the Corollary, in particular,
with (some of) the equilibria replaced by periodic orbits.
We intend to address this question in the future.

For heteroclinic cycles the stability indices provide quantitative and
qualitative description for behaviour of nearby trajectories. In a f.a.s.
heteroclinic network the trajectory through a point that belongs to its local
basin can be possibly
attracted by any of its f.a.s. subcycles, or it can switch between different
subcycles, without being attracted by any of them. Certainly,
stability indices, either for the whole network or for individual cycles,
do not provide such information. One can think about proposing for
a network $X$, a subset $Y\subset X$ and a point $x\in X$ a relative stability
index $\sigma(x,Y,X)$, describing the part of a small neighbourhood of
$x$ that stays near $X$ for all $t>0$ and is attracted by $Y$ as $t\to\infty$.
If in addition the sets $Y$ are required to be maximal and undecomposable,
such relative stability indices should be useful for describing local dynamics
near the network.
Introduction of such an index is beyond the scope of the present paper.
A step toward this direction
was made in \cite{cl2014} by defining stability indices with respect to a
cycle (the $c$-index) and with respect to the whole network (the $n$-index).

\section*{Acknowledgements}

All  authors were partially supported by CMUP (UID/MAT/00144/2013), which is funded by FCT (Portugal) with national (MEC) and European structural funds (FEDER), under the partnership agreement PT2020.
Much of the work was done while O.P. was visiting CMUP, whose hospitality is gratefully acknowledged.

\appendix

\section{Lemmas in the proof of Theorem~\ref{th1}}\label{appB}

Denote by $\tilde{g}_j^{(k)}$ , where $j=1,2,3,4$,
the $j$-th component of the map $\psi_{ki}\phi_k$, with $i=k+1\pmod{3}$.

\begin{lemma}\label{lem51}
For any $\delta>0$ there exist $\varepsilon>0$ such that
$|{u}_{35}^\ini|<\varepsilon$ implies that
$$
(i)\qquad
|A_4\tilde g_3^{(3)}({u}_{31}^\ini,{u}_{32}^\ini,{u}_{33}^\ini,{u}_{35}^\ini)|-
|A_3\tilde g_4^{(3)}({u}_{31}^\ini,{u}_{32}^\ini,{u}_{33}^\ini,{u}_{35}^\ini)|<
\delta|\tilde g_3^{(3)}({u}_{31}^\ini,{u}_{32}^\ini,{u}_{33}^\ini,{u}_{35}^\ini)|
$$
and
$$
(ii)\qquad
|\tilde g_3^{(3)}({u}_{31}^\ini,{u}_{32}^\ini,{u}_{33}^\ini,{u}_{35}^\ini)-
\tilde g_3^{(3)}({u}_{31}^\ini,{u}_{32}^\ini,{u}_{33}^\ini,0)|<
\delta|\tilde g_3^{(3)}({u}_{31}^\ini,{u}_{32}^\ini,{u}_{33}^\ini,0)|.
$$
\end{lemma}

\proof
 Recall that $\lambda_{35}=\lambda_{36}$.
Denote $(x_1,x_2,x_3,x_4)=({u}_{32}^\ini,{u}_{33}^\ini,{u}_{31}^\ini,{u}_{35}^\ini)$.
From (\ref{psi})--(\ref{phi}) and Table~\ref{tbBases} we have
$$
\tilde g_3^{(3)}(x_1,x_2,x_3,x_4)=
 \frac{A_3}{\sqrt{2}}
(x_4|x_3|^{-\lambda_{35}/\lambda_{31}}+
D_1|x_3|^{-\lambda_{35}/\lambda_{31}})
$$
and
$$
\tilde g_4^{(3)}(x_1,x_2,x_3,x_4)=
\frac{A_4}{\sqrt{2}}
(x_4|x_3|^{-\lambda_{35}/\lambda_{31}}-
D_1|x_3|^{-\lambda_{35}/\lambda_{31}}).
$$

To prove (i) note that $|A_4\tilde g_3^{(3)}|-|A_3\tilde g_4^{(3)}|\le |A_4\tilde g_3^{(3)}+A_3\tilde g_4^{(3)}|$.
Then
$$
|A_4\tilde g_3^{(3)}(x_1,x_2,x_3,x_4)+A_3\tilde g_4^{(3)}(x_1,x_2,x_3,x_4)|=
2A_4 \frac{A_3}{\sqrt{2}}|x_4||x_3|^{-\lambda_{35}/\lambda_{31}}
$$
and
$$
\delta | \tilde g_3^{(3)}(x_1,x_2,x_3,x_4)|= \delta \frac{A_3}{\sqrt{2}} |x_4+D_1| |x_3|^{-\lambda_{35}/\lambda_{31}}.
$$
Note that
$$
2A_4|x_4| < \delta |D_1|-\delta |x_4|\le \delta |D_1+x_4|
$$
and the first inequality holds if $\displaystyle |x_4| < \frac{\delta D_1}{2A_4+\delta}$.
Substituting in (ii) we obtain $|x_4|< \delta |D_1|$.
To finish the proof choose $\displaystyle \varepsilon= \min\{ \frac{\delta D_1}{2A_4+\delta}, \delta |D_1|\}$.
\qed

Lemma~\ref{lem51} implies that for trajectories sufficiently close to
the cycle the value of $u_{35}^\ini$ is irrelevant. Therefore,
in the study of stability, for simplicity,
we can ignore $u_{35}^\ini$ and instead of
$\tilde g:\,\R^4\to\R^4$ (which maps $H^\ini_3\to H^\ini_3$) consider $\bg:\R^3\to\R^3$,
\begin{equation}\label{defg}
\begin{array}{lcl}
\bg(x_1,x_2,x_3)&=&(\tilde g_1(x_1,x_2,x_3,0),
\tilde g_2(x_1,x_2,x_3,0),\tilde g_3(x_1,x_2,x_3,0))\\
&=&\left(g_1(x_1,x_2,x_3),g_2(x_1,x_2,x_3),g_3(x_1,x_2,x_3) \right).
\end{array}
\end{equation}
From (\ref{phi})--(\ref{psi}) and Lemma~\ref{lem51}
choosing coordinates in each $H_j^\ini$:
\begin{equation}\label{varrr}
\renewcommand{\arraystretch}{1.5}
\begin{array}{cl}
H_3^{\ini}:& (x_1,x_2,x_3)=(u_{32}^\ini,u_{33}^\ini,u_{31}^\ini)\\
H_1^{\ini}:& (x_1,x_2,x_3)=(u_{12}^\ini,u_{13}^\ini,u_{16}^\ini=-A_3u_{15}^\ini/A_4)\\
H_2^{\ini}:& (x_1,x_2,x_3)=(u_{22}^\ini,u_{23}^\ini,u_{26}^\ini)
\end{array}
\end{equation}
and using $\lambda_{12}=\lambda_{13}$, $\lambda_{22}=\lambda_{23}$ and $\lambda_{32}=\lambda_{33}$,
we can write
$\bg=\bg^{(2)}\bg^{(1)}\bg^{(3)}$, where
\begin{equation}\label{mapg}
\renewcommand{\arraystretch}{1.5}
\begin{array}{cl}
\bg^{(3)}(x_1,x_2,x_3)=&\frac{\sqrt{2}}{2}
\left(A_1(x_1+x_2)|x_3|^{-\lambda_{32}/\lambda_{31}},\right.\
A_2(x_1-x_2)|x_3|^{-\lambda_{32}/\lambda_{31}},\\
&\left.A_3D_1|x_3|^{-\lambda_{35}/\lambda_{31}}\right)\\
\bg^{(1)}(x_1,x_2,x_3)=&\left(B_1x_1|A_4x_3/A_3|^{-\lambda_{12}/\lambda_{15}},\right.\
B_2x_2|A_4x_3/A_3|^{-\lambda_{12}/\lambda_{15}},\\
&\left.-B_4x_3|A_4x_3/A_3|^{-\lambda_{16}/\lambda_{15}}\right)\\
\bg^{(2)}(x_1,x_2,x_3)=&\left(C_3{D_3}|x_3|^{-\lambda_{21}/\lambda_{26}},\
(C_4x_1+C_5x_2)|x_3|^{-\lambda_{22}/\lambda_{26}},\right.\\
&\left.(C_1x_1+C_2x_2)|x_3|^{-\lambda_{22}/\lambda_{26}}\right).
\end{array}
\end{equation}

From (\ref{phi}) and (\ref{psi}) for $\psi_{23}\phi_2$ we have:
$$
\renewcommand{\arraystretch}{1.5}
\begin{array}{l}
\hat u_{31}^\ini=C_1u_{22}^\ini|u_{26}^\ini|^{-\lambda_{22}/\lambda_{26}}+
C_2u_{23}^\ini|u_{26}^\ini|^{-\lambda_{22}/\lambda_{26}},\\
\hat u_{32}^\ini=C_3D_3|u_{26}^\ini|^{-\lambda_{21}/\lambda_{26}},\\
\hat u_{33}^\ini=C_4u_{22}^\ini|u_{26}^\ini|^{-\lambda_{22}/\lambda_{26}}+
C_5u_{23}^\ini|u_{26}^\ini|^{-\lambda_{22}/\lambda_{26}}.
\end{array}
$$
Then we recall that $(x_1,x_2,x_3)=(u_{22}^\ini,u_{23}^\ini,u_{26}^\ini)$
and $(\hat u_{32}^\ini,\hat u_{33}^\ini,\hat u_{31}^\ini)=(x_1,x_2,x_3)$, see (\ref{varrr}).

Therefore, we obtain that
\begin{equation}\label{mapg0}
\renewcommand{\arraystretch}{1.5}
\begin{array}{cl}
g_1(x_1,x_2,x_3)=&F_1|x_3|^{\beta_1},\\
g_2(x_1,x_2,x_3)=&(F_2x_1+F_3x_2)|x_3|^{\beta_2},\\
g_3(x_1,x_2,x_3)=&(F_4x_1+F_5x_2)|x_3|^{\beta_2},\\
\end{array}
\end{equation}
where
\begin{equation}\label{ab1122}
\beta_1=\frac{\lambda_{21}\lambda_{35}}{\lambda_{26}\lambda_{31}}(1-\frac{\lambda_{16}}{\lambda_{15}}),\qquad
\beta_2=-\frac{\lambda_{32}}{\lambda_{31}}+
\frac{\lambda_{35}\lambda_{12}}{\lambda_{31}\lambda_{15}}+
\frac{\lambda_{35}\lambda_{22}}{\lambda_{31}\lambda_{26}} \left(1-\frac{\lambda_{16}}{\lambda_{15}}\right)
\end{equation}
and $F_j$ depend on $A_i$, $B_i$, $C_i$ and $D_i$. Generically, $F_j\ne0$,
$F_2\ne F_3$, $F_4\ne F_5$ and $F_2/F_3\ne F_4/F_5$. For definiteness
we assume that all $F_j$ are positive.

 Evidently, $\beta_1<0$ implies that the origin is a
completely unstable fixed point of the map $\bg$ (\ref{mapg0}).
From now on till the end of this subsection we assume that $\beta_1>0$.

Given the map $\bg$ (\ref{mapg0}), we define the map $\bh:~\R_+^2\to\R_+^2$ as
$$
\bh(p,q)=(\max\{pq^{\beta_2},q^{\beta_1}\},pq^{\beta_2}).
$$
This corresponds to the map  in Lemmas~\ref{lem4new} and \ref{lem4}, with parameters:
$$
\alpha_1=0\quad \alpha_2=1\quad\gamma=1\quad \gamma_1=\beta_1-\beta_2.
$$
In Lemmas~\ref{lem7}--\ref{lem84}
we prove that the stability properties of the
origin, which is a fixed point of both $\bg$ and $\bh$, are
the same for both maps, namely that the origin is
either a.s. if $\beta_1>0$, $\beta_2>0$ and $\beta_1+\beta_2>1$, or it is
c.u. otherwise (see Lemmas~\ref{lem4new} and \ref{lem4}).
A hand-waving proof can be obtained by denoting
\begin{equation}\label{nott}
p=\max\{|x_1|,|x_2|\}
\quad\hbox{ and }\quad
q=|x_3|,
\end{equation}
ignoring constants in (\ref{mapg0})
and noting that for small $\bf x$ we have generically either
$x_3^{\beta_1}\gg\max\{|x_1|,|x_2|\}x_3^{\beta_2}$ or
$x_3^{\beta_1}\ll\max\{|x_1|,|x_2|\}x_3^{\beta_2}$.
A rigorous proof is given below in a series of lemmas.

\begin{lemma}\label{lem7}
If $\beta_1>0$, $\beta_2>0$ and $\beta_1+\beta_2>1$ then the origin
is an asymptotically stable fixed point of the map $\bg$ (\ref{mapg0}).
\end{lemma}

\proof
There exist $s>0$ and $\varepsilon>0$ such that $\beta_1-s>0$, $\beta_2-s>0$,
$\beta_1+\beta_2-2s>1$ and $\max\{|F_1|,|F_2|+|F_3|,|F_4|+|F_5|\}\varepsilon^s<1$.
According to Lemmas~\ref{lem4new} and \ref{lem4}, the origin is an asymptotically
stable fixed point of the map
$$\bh^s(p,q)=(\max\{pq^{\beta_2-s},q^{\beta_1-s}\},pq^{\beta_2-s}).$$
For $|\bx|<\varepsilon$ we have
$$|g_1(x_1,x_2,x_3)|<q^{\beta_1-s},\
|g_2(x_1,x_2,x_3)|<pq^{\beta_2-s},\
|g_3(x_1,x_2,x_3)|<pq^{\beta_2-s},$$
where $p$ and $q$ are defined by (\ref{nott}). Since
$p_0\le p$ and $q_0\le q$ imply that $h_1^s(p_0,q_0)\le h_1^s(p,q)$ and
$h_2^s(p_0,q_0)\le h_2^s(p,q)$, for any $n>0$ and $|\bx|<\varepsilon$
the iterates $\bg^n(\bx)$ satisfy
$$|g_1^n(x_1,x_2,x_3)|<(h^s)_1^n(p,q),\ |g_2^n(x_1,x_2,x_3)|<(h^s)_1^n(p,q),\
|g_3^n(x_1,x_2,x_3)|<(h^s)_2^n(p,q).$$
\qed

\begin{lemma}\label{lem81}
Consider the map $\bg$ (\ref{mapg0}) where $\beta_1>0$ and $\beta_2<-1$.
The origin is completely unstable fixed point of the map $\bg$ (\ref{mapg0}).
\end{lemma}

\begin{figure}
\centerline{
\includegraphics[scale=1.1]{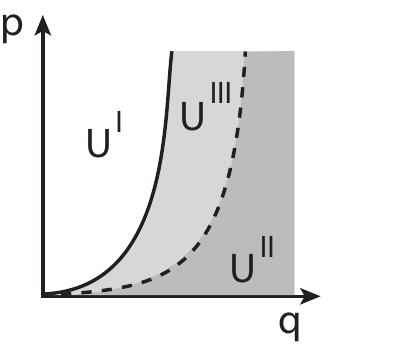}
}
\caption{
Decomposition for Lemma~\ref{lem81} in the $(p,q)$ plane. The dashed line,  the common boundary of $U^{II}$ and  $U^{III}$, is mapped by $h^{II}$ into the solid line,  the common boundary of $U^{I}$ and  $U^{III}$.
\label{setUlemB3}}
\end{figure}

\proof
 Let
$
\R^3_{\delta}=\{(x_1,x_2,x_3)~:~\max\{p,q\}<{\delta}\} .
$
We decompose
$\R^3_{\delta}=U^I\cup U^{II}\cup U^{III}$ and
$U^{III}=U^{III}_{c_0}\cup W_{c_0}$, where
\begin{equation}\label{defsets}
\renewcommand{\arraystretch}{1.5}
\begin{array}{l}
U^I=\{(x_1,x_2,x_3)~:~\max\{p,q\}<{\delta}\hbox{ and }p>q^{-\beta_2}\},\\
U^{II}=\{(x_1,x_2,x_3)~:~\max\{p,q\}<{\delta}\hbox{ and }
p<q^{-(\beta_1+\beta_2^2)/\beta_2}\},\\
U^{III}=\{(x_1,x_2,x_3)~:~\max\{p,q\}<{\delta}\hbox{ and }
q^{-(\beta_1+\beta_2^2)/\beta_2}<p<q^{-\beta_2}\},\\
U^{III}_{c_0}=\{(x_1,x_2,x_3)\in U^{III}~:~|F_2x_1+F_3x_2|>c_0|x_1|\},\
W_{c_0}=U^{III}\setminus U^{III}_{c_0},
\end{array}
\end{equation}
$0<c_0<\min\{1,F_2\}/2$ and as above $p=\max\{|x_1|,|x_2|\}$ and $q=|x_3|$.

As before, let
$\bh^I(p,q)=(pq^{\beta_2},pq^{\beta_2})$ and
$\bh^{II}(p,q)=(q^{\beta_1},pq^{\beta_2})$.
The curve $p=q^{-(\beta_1+\beta_2^2)/\beta_2}$,
the common boundary of $U^{II}$ and  $U^{III}$, is mapped by $h^{II}$ into the common boundary of $U^{I}$ and  $U^{III}$, the curve $p=q^{-\beta_2}$.

For sufficiently small $\delta$ we claim that
the subsets are mapped by $\bg$ as follows:
\begin{equation}\label{mapsets}
\renewcommand{\arraystretch}{1.5}
\bg (U^I)\cap V_{\delta}=\varnothing;\qquad
\bg (U^{II})\cap V_{\delta}\subset U^I;\qquad
\bg (U^{III}_{c_0})\cap V_{\delta}\subset U^I.
\end{equation}

The first equality holds  under the generic assumption $F_2/F_3\ne F_4/F_5$.
To see this, let $\dpt C=\min_{p=1} \max\{ |F_2x_1+F_3x_2|, |F_4x_1+F_5x_2|\}>0$.
The genericity assumption guarantees that $C\ne 0$.
Then for  $p<\delta$ we have
$$
 \max\{|F_2x_1+F_3x_2|, |F_4x_1+F_5x_2|\}> Cp.
 $$
When $\bx\in U^I$, since $p>q^{-\beta_2} $ this implies
either $|\bg_2(\bx)|>C$ or $|\bg_3(\bx)|>C$ and  the assertion holds.

For the image of $U^{II}$, notice that the map $\bg_3$ is a decreasing function of $x_3$, i.e. if $x_3< \hat{x}_3$, then $\bg_3(x_1,x_2,\hat{x}_3)<\bg_3(x_1,x_2,x_3)$.

The last inclusion holds because for $\bx\in U^{III}_{c_0}$
we have $c_0|x_1|> \left| F_3x_2+F_2x_1\right|\ge F_3|x_2|-F_2|x_2|$  and hence
 $\displaystyle |x_1|>\frac{F_3}{F_2+c_0}|x_2|$, therefore
$$
|g_2(\bx)|>\min\{c_0,\frac{c_0F_3}{ F_2+c_0}\} pq^{\beta_2}
\quad\hbox{ and }\quad
|g_3(\bx)|<|F_4+F_5|pq^{\beta_2}.
$$

\begin{figure}
\centerline{
\includegraphics[scale=1.2]{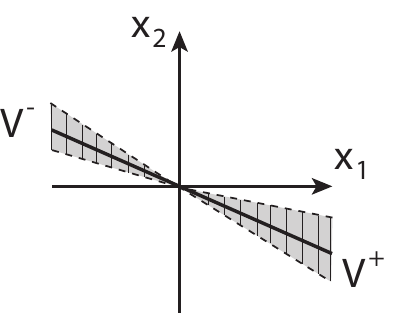}
\qquad\qquad
\includegraphics[scale=1.1]{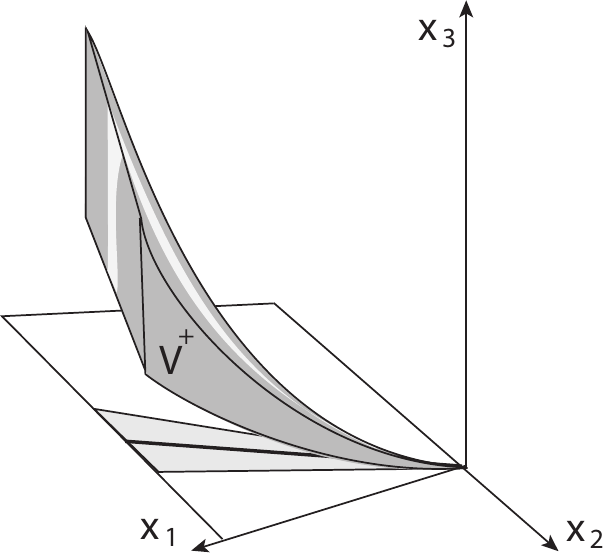}
}
\caption{
On the left:  Projection into the $(x_1,x_2)$-plane of the two sets $V^\pm$ that are bissected by the solid line $F_2x_1+F_3 x_2=0$ and are covered by the lines $x_1=$constant, shown here in the case $0<F_2/F_3<1$. Right: the set $V^+$ and its projection. Only the lighter strip is mapped back into $V^+$ by $\bg$.
\label{setW+}}
\end{figure}

To investigate how $\bx\in W_{c_0}$ are mapped by $\bg$, we decompose
$W_{c_0}=W^+\cup W^-$ where $W^i=\{(x_1,x_2,x_3)\in W_{c_0}: ~  i x_1>0\}$ and
$$
W^+=
\bigcup_{\stackrel{  -\beta_2/(\beta_1+\beta_2^2) <r< -1/\beta_2}{x_1<\delta}   }
V^{+}(x_1,- c_0, c_0,r),
\qquad
W^-=
\bigcup_{\stackrel{  -\beta_2/(\beta_1+\beta_2^2) <r< -1/\beta_2}{-x_1<\delta}   }
V^{-}(x_1,c_0, -c_0,r),
$$
where
\begin{equation}\label{defV}
\begin{array}{l}
V^+(x_1,c_1,c_2,r)=\{(x_1,x_2,x_1^r)~:~c_1x_1<|F_2x_1+F_3x_2|<c_2x_1,~x_1>0\},\\
V^-(x_1,c_1,c_2,r)=\{(x_1,x_2,x_1^r)~:~c_1x_1<|F_2x_1+F_3x_2|<c_2x_1,~x_1<0\}.
\end{array}
\end{equation}
We have:
$$
\bg(x_1,x_2,x_1^r)=
(F_1x_1^{r\beta_1},(F_2x_1+F_3x_2)x_1^{r\beta_2},(F_4x_1+F_5x_2)x_1^{r\beta_2}).
$$
 For fixed $x_1$ and $r$, look at the curve $x_2\mapsto(x_1,x_2, x_1^r)$ inside $W^+$.
We claim that only a small part of this curve is mapped by $\bg$ into $W_{c_0}$.
The set $\bg V^+(x_1,-c_0,c_0,r)\cap W_{c_0}$ satisfies
$$-c_0F_1x_1^{r\beta_1}<
F_2F_1x_1^{r\beta_1}+
F_3(F_2x_1+F_3x_2)x_1^{r\beta_2}
<c_0F_1x_1^{r\beta_1},$$
which implies that
the points that remain in $W_{c_0}$ are
\begin{equation}\label{gVpl}
\bg^{-1}(\bg (V^+(x_1,-c_0,c_0,r))\cap W_{c_0})
=V^+(x_1,c'-\tilde c,c'+\tilde c,r),
\end{equation}
where
\begin{equation}\label{ctld}
c'=-\frac{F_1F_2}{ F_3}x_1^{\beta_1r-\beta_2r-1},\quad
\tilde c=-\frac{c_0F_1}{F_3}x_1^{\beta_1r-\beta_2r-1}.
\end{equation}
(Note, that $r>-\beta_2/(\beta_1+\beta_2^2)$ implies that
$\beta_1r-\beta_2r-1>\beta_1(-\beta_2-1)/(\beta_1+\beta_2^2)>0$.)
Moreover, for asymptotically small $x_1$ we can write
\begin{equation}\label{mapVp}
\bg (V^+(x_1,-c_0,c_0,r))\cap W_{c_0}\approx
V^+(F_1x_1^{r\beta_1},-c_0,c_0,r'),\hbox{ where } r'=(1+\beta_2r)/(r\beta_1).
\end{equation}
 Here $r^\prime>0$ because $r<-1/\beta_2$.
Similar estimates holds true for $V^-(x_1,c_0,-c_0,r)$.

 We represent each one of the sets $W^+$ and $W^-$ as the union of a $(x_1, r)$-familiy of curves parametrised by $x_2$.
From (\ref{gVpl}) and (\ref{ctld}) we obtain that
\begin{equation}\label{estw1}
\frac{l_1(\bg^{-1}(\bg V^+(x_1,-c_0,c_0,r)\cap W_{c_0}))}{
l_1(V^+(x_1,-c_0,c_0,r)))}=\frac{\tilde c}{ c_0}<Cx_1^{\beta_1r-\beta_2r-1},
\end{equation}
where $l_1$ denotes the 1-dimensional Lebesgue measure and $C$ is a
constant that depends on $F_j$ and $\beta_i$. Hence, for
sufficiently small $\delta$ the inclusions (\ref{mapsets}) and
the estimate (\ref{mapVp}) imply that as $n\to\infty$ almost all points, except
for a set of zero measure, are mapped by $\bg^n$ away from $V_{\delta}$.
Therefore, the point $\bx=0$ is a completely unstable point
of the map $\bg$.
\qed

\begin{lemma}\label{lem82}
Consider the map $\bg$ (\ref{mapg0}) where $\beta_1>0$, $\beta_1-\beta_2>1$
and $-1<\beta_2<0$.
The origin is completely unstable fixed point of the map $\bg$ (\ref{mapg0}).
\end{lemma}

\proof
The proof of this lemma, and also of the two following, employs the same ideas
as the proof of Lemma~\ref{lem81}. Namely, we decompose $\R^3_{\delta}$ as
a union of several subsets and consider how the subsets are mapped by $\bg$.
We represent $\R^3_{\delta}=U^I\cup U^{II}$ and $U^I=U^I_{c_0}\cup W_{c_0}$, where
\begin{equation}\label{defsets2}
\renewcommand{\arraystretch}{1.5}
\begin{array}{l}
U^I=\{(x_1,x_2,x_3)~:~\max\{p,q\}<{\delta}\hbox{ and }p>q^{\beta_1-\beta_2-s}\},\\
U^{II}=\{(x_1,x_2,x_3)~:~\max\{p,q\}<{\delta}\hbox{ and }p<q^{\beta_1-\beta_2-s}\},\\
U^I_{c_0}=\{(x_1,x_2,x_3)\in U^I~:~|F_2x_1+F_3x_2|>c_0|x_1|\},\
W_{c_0}=U^I\setminus U^I_{c_0},
\end{array}
\end{equation}
$0<s<\beta_1$, $\beta_1-\beta_2-s>\beta_1/(\beta_1-s)$ and
$0<c_0<\min\{1,F_2\}/2$.

For sufficiently small $\delta$ we have:
\begin{equation}\label{mapsets2}
(\bg( U^I_{c_0} )\cap V_{\delta} )\subset U^I_{c_0};\qquad
(\bg (U^{II} )\cap V_{\delta} )\subset U^I.
\end{equation}
The latter inclusion holds true due to our choice of $s$.

For any $\bx\in U^I_{c_0}$ for small $\delta$ we have
$g_1(\bx)\ll g_2(\bx)$, hence for $n>2$ we can write
$$\bg^n(\bx)\approx
(F_1 (\bg^{n-1}_3 (\bx))^{\beta_1},
F_3 \bg^{n-1}_2(\bg^{n-1}_3 (\bx))^{\beta_2},
F_5 \bg^{n-1}_2(\bg^{n-1}_3 (\bx))^{\beta_2}).$$
Therefore,
$$\bg^{n+1}(\bx)\approx
(F_1 (\bg^n_3 (\bx))^{\beta_1},
F_3' (\bg^n_3 (\bx))^{1+\beta_2},
F_5' (\bg^n_3 (\bx))^{1+\beta_2}).$$
which implies that for any $\bx\in U^I_{c_0}$ we have
$\bg^n(\bx)\not\in V_{\delta}$ for large $n$.

We decompose further:
\begin{equation}\label{decow}
W_{c_0}=\cup_{\pm}
\bigcup_{x_1<\delta,~-1/(\beta_1-\beta_2-s)<r<\infty} V^{\pm}(x_1,\mp c_0,\pm c_0,r),
\end{equation}
where $V^{\pm}(x_1,c_1,c_2,r)$ are defined by (\ref{defV}).

The set $V^+$ is mapped as $\bg(V^+)\subset U^I_{c_0}\cup W_{c_0}\cup U^{II}$.
We have accounted for points that go to $U^I_{c_0}$.
Below, we first show that the set of points that are mapped to $W_{c_0}$ is small.
Then, for points mapped to $U^{II}$ we will have to consider a second iteration of $\bg$
to show that the set of points that are first mapped to  $U^{II}$  and then to $W_{c_0}$ is small.
(The points that are first mapped to $U^{II}$ and then to $U^I_{c_0}$ are already accounted for. )
As in the proof of the previous lemma, this implies that when the number of iterations goes to infinity
almost all points are mapped to  $U^I_{c_0}$, and then away from $V_\delta$.

We have (see the proof of Lemma~\ref{lem81})
\begin{equation}\label{gVpl2}
\bg^{-1}(\bg( V^+(x_1,-c_0,c_0,r) )\cap W_{c_0})
=V^+(x_1,c'-\tilde c,c'+\tilde c,r),
\end{equation}
where $\tilde c=c_0F_1x_1^{\beta_1r-\beta_2r-1}/F_3$.
Since $\beta_1r-\beta_2r-1>s/(\beta_1-\beta_2-s)>0$, we estimate
\begin{equation}\label{estw2}
\frac{l_1(\bg^{-1}(\bg (V^+(x_1,-c_0,c_0,r) )\cap W_{c_0}))}{
l_1(V^+(x_1,-c_0,c_0,r)))}=\frac{\tilde c}{ c_0}< Cx_1^{s/(\beta_1-\beta_2-s)}.
\end{equation}

Next, we consider $\bg(\bg( V^+(x_1,-c_0,c_0,r) )\cap U^{II})$.
For $\bx=(x_1,-F_2x_1/F_3+cx_1,x_1^r)$ we have
$$\bg^2(\bx)=
(\tilde F_1x_1^{\beta_1(1+r\beta_2)},
(\tilde F_2x_1^{\beta_1r}+c\tilde F_3x_1^{1+r\beta_2})x_1^{\beta_2(1+r\beta_2)},
(\tilde F_4x_1^{\beta_1r}+c\tilde F_5x_1^{1+r\beta_2})x_1^{\beta_2(1+r\beta_2)}).$$
Therefore,
\begin{equation}\label{estww2}
\bg^{-2}(\bg(\bg( V^+(x_1,-c_0,c_0,r) )\cap U^{II})\cap W_{c_0})=
V^+(x_1,c^*-\hat c,c^*-\hat c,r),
\end{equation}
where $\hat c=\hat C x_1^{(1+r\beta_2)(\beta_1-\beta_2-1)}$.
Moreover,
\begin{equation}\label{estw22}
\bg(\bg (V^+(x_1,-c_0,c_0,r) )\cap U^{II})\cap W_{c_0}\subset
V^+(x_1,-c_0,c_0,r'),
\end{equation}
where $r'=(\min\{r\beta_1,1+r\beta_2\}+\beta_2)/\beta_1$.
Therefore, due to (\ref{mapsets2})-(\ref{estw22}) for sufficiently
small $\delta$, almost all $\bx\in V_{\delta}$ satisfy
$|\bg^n(\bx)|>\delta$ as $n\to\infty$.
\qed

\begin{lemma}\label{lem83}
Consider the map $\bg$ (\ref{mapg0}) where $\beta_1>0$, $\beta_1-\beta_2<1$
and $-1<\beta_2<0$.
The origin is completely unstable fixed point of the map $\bg$ (\ref{mapg0}).
\end{lemma}

\proof
We decompose
$\R^3_{\delta}=U^I\cup U^{II}\cup U^{III}$ and $U^I=U^I_{c_0}\cup W_{c_0}$, where
\begin{equation}\label{defsets3}
\renewcommand{\arraystretch}{1.5}
\begin{array}{l}
U^I=\{(x_1,x_2,x_3)~:~\max\{p,q\}<\delta\hbox{ and }p>q^{\beta_1-\beta_2-s_1}\},\\
U^{II}=\{(x_1,x_2,x_3)~:~\max\{p,q\}<\delta\hbox{ and }p<q^{1-s_2}\},\\
U^{III}=\{(x_1,x_2,x_3)~:~\max\{p,q\}<\delta\hbox{ and }
q^{1-s_2}<p<q^{\beta_1-\beta_2-s_1}\},\\
U^I_{c_0}=\{(x_1,x_2,x_3)\in U^I~:~|F_4x_1+F_5x_2|>c_0|x_1|\},\
W_{c_0}=U^{I}\setminus U^{I}_{c_0},
\end{array}
\end{equation}
$0<s_1<
\min\{1,\beta_1\beta_2(\beta_1-\beta_2-1)/(\beta_1-\beta_2(\beta_1-\beta_2))\}$,
$0<s_2<\min\{1,\beta_2(\beta_1-\beta_2-1)/(\beta_1-\beta_2)\}$
and $c_0<\min\{1,F_4\}/2$.

Due to our choice of $s_1$, $s_2$ and $c_0$, for
sufficiently small $\delta$ the subsets are mapped by $\bg$ as follows:
\begin{equation}\label{mapsets3}
\left(\bg( U^I_{c_0} )\cap V_{\delta } \right)\subset U^{II};\qquad
\left(\bg( U^{II} )\cap V_{\delta} \right)\subset U^I_{c_0}.
\end{equation}

Consider $\bx\in U^I_{c_0}\cup U^{II}$. Due to (\ref{mapsets3}),
we can assume that $\bx\in U^{II}$. Denote
$p_n=\max\{|(\bg^n(\bx))_1|,|(\bg^n(\bx))_2|\}$ and $q_n=|(\bg^n(\bx))_3|$.
For $n\to\infty$ the iterates $\bg^n(\bx)$ as long as
$\bg^n(\bx)\in V_{\delta}$ satisfy\\
for even $n$: $(p_n,q_n)\approx(F_1'q_{n-1}^{\beta_1},F_2'p_{n-1}q_{n-1}^{\beta_1})$;\\
for odd $n$:
$(p_n,q_n)\approx(F_3'p_{n-1}q_{n-1}^{\beta_1},F_4'p_{n-1}q_{n-1}^{\beta_1})$.\\
Which implies that
$$p_n\approx F_5'q_n
\quad\hbox{ and }\quad
q_n\approx F_6'q_{n-2}^{\beta_1+\beta_2(1+\beta_2)}.$$
By assumption, $\beta_1>0$ and $\beta_1-\beta_2<1$, which implies that
$\beta_1+\beta_2(1+\beta_2)<1$, hence
the iterates $\bg^n(\bx)$ satisfy $|\bg^n(\bx)|>\delta$ as $n\to\infty$.

On the other hand, for sufficiently small $\bx\in U^{III}$ the map $\bg$ can
be approximated by $\bh(p,q)=(F_1'q^{\beta_1},F_2'pq^{\beta_2})$. Therefore,
arguments similar to the ones employed in the proof of Lemma~\ref{lem4}
imply that for almost all $\bx$ the iterates $\bg^n(\bx)$ escape
from $U^{III}$ as $n\to\infty$. The iterates satisfy:
\begin{equation}\label{mapsets2a}
\hbox{if }\bx\in U^{III}\hbox{ and }\bg(\bx)\in U^{III}\hbox{ then }
\bg^2(\bx)\not\in W_{c_0}.
\end{equation}
By decomposing $W_{c_0}$ similarly to (\ref{decow}), proceeding as in
Lemma~\ref{lem82} by considering $\bg(\bx)$ and $\bg^2(\bx)$ for $\bx\in W_{c_0}$
and taking into account the above inclusions, we obtain that almost
all points in $V_{\delta}$ either  escape from $V_{\delta}$ or are
mapped by $\bg^2$ to $U^I_{c_0}\cup U^{II}$, from where they escape from $V_{\delta}$.
\qed

\begin{lemma}\label{lem84}
Consider the map $\bg$ (\ref{mapg0}) where $\beta_1>0$, $\beta_2>0$ and
$\beta_1+\beta_2<1$.
The origin is completely unstable fixed point of the map $\bg$ (\ref{mapg0}).
\end{lemma}

\proof
We decompose
$\R^3_{\delta}=U^I\cup U^{II}\cup U^{III}$, $U^I_{c_0}\subset U^I$,
$U^{III}_{c_0}\subset U^{III}$ and
$W_{c_0}=(U^I\cup U^{III})\setminus (U^I_{c_0}\cup U^{III}_{c_0})$, where
\begin{equation}\label{defsets4}
\renewcommand{\arraystretch}{1.5}
\begin{array}{l}
U^I=\{(x_1,x_2,x_3)~:~\max\{p,q\}<\delta\hbox{ and }p>q^{\beta_1-\beta_2-s_1}\},\\
U^{II}=\{(x_1,x_2,x_3)~:~\max\{p,q\}<\delta\hbox{ and }p<q^{1+s_2}\},\\
U^{III}=\{(x_1,x_2,x_3)~:~\max\{p,q\}<\delta\hbox{ and }
q^{1+s_2}<p<q^{\beta_1-\beta_2-s_1}\},\\
U^I_{c_0}=\{(x_1,x_2,x_3)\in U^I~:~|F_4x_1+F_5x_2|>c_0|x_1|\},\\
U^{III}_{c_0}=\{(x_1,x_2,x_3)\in U^{III}~:~|F_4x_1+F_5x_2|>c_0|x_1|\},
\end{array}
\end{equation}
with $0<s_1<
\min\{1,\beta_2(1-\beta_1+\beta_2)/(1+\beta_2))\}$,
$0<s_2<\min\{1,\beta_2(1-\beta_1+\beta_2)/(\beta_1-\beta_2)\}$
and $c_0<\min\{1,F_4\}/2$.

For sufficiently small $\delta$ we have:
\begin{equation}\label{mapsets4}
\left( \bg( U^{III}_{c_0}\ )cap V_{\delta} \right) \subset U^{III}_{c_0};\qquad
\left( \bg (U^I_{c_0} )\cap V_{\delta} \right) \subset U^{III}_{c_0};\qquad
\left( \bg (U^{II} )\cap V_{\delta} \right) \subset U^I_{c_0}.
\end{equation}
By decomposing $W_{c_0}$ similarly to (\ref{decow}) and considering
$\bg(\bx)$ and $\bg^2(\bx$) for $\bx\in W_{c_0}$, we prove that
all points in $V_{\delta}$ are either mapped by $\bg^2$ to $U^{III}_{c_0}$,
or they escape from $V_{\delta}$.
For $\bx\in U^{III}_{c_0}$ we note that the map $\bg(\bx)$ can be approximated
by $\bh(p,q)=(F_1'q^{\beta_1},F_2'pq^{\beta_2})$, which implies that
the iterates $\bg^n(\bx)$ satisfy $|\bg^n(\bx)|>\delta$ as $n\to\infty$.
\qed

\section{Eigenspaces and eigenvalues
near single-mode steady states}\label{AppTable}

In this appendix we provide data on the network that are used in the calculations.

\begin{table}
$$
\renewcommand{\arraystretch}{1.5}
\begin{array}{l}
\be_{11}=(1,0,0;0,0,0)/\sqrt{2},\ \be_{12}=(0,1,1;0,0,0)/\sqrt{2},\
\be_{13}=(0,1,-1;0,0,0)/\sqrt{2},\\
\be_{14}=(0,0,0;1,0,0),\
\be_{15}=(0,0,0;0,1,1)/\sqrt{2},\ \be_{16}=(0,0,0;0,1,-1)/\sqrt{2}\\
\hline
\be_{21}=(1,0,0;0,0,0),\ \be_{22}=(0,1,1;0,0,0)/\sqrt{2},
\ \be_{23}=(0,1,-1;0,0,0)/\sqrt{2}\\
\be_{24}=(0,0,0;1,0,0),\ \be_{25}=(0,0,0;0,1,1)/\sqrt{2},\
\be_{26}=(0,0,0;0,1,-1)/\sqrt{2}\\
\hline
\be_{31}=(1,0,0;0,0,0),\ \be_{32}=(0,1,0;0,0,0),\ \be_{33}=(0,0,1;0,0,0),\\
\be_{34}=(0,0,0;1,0,0),\ \be_{35}=(0,0,0;0,1,0),\ \be_{36}=(0,0,0;0,0,1)\\
\tilde\be_{31}=(1,0,0;0,0,0),\ \tilde\be_{32}=(0,1,1;0,0,0)/\sqrt{2},
\ \tilde\be_{33}=(0,1,-1;0,0,0)/\sqrt{2}\\
\tilde\be_{34}=(0,0,0;1,0,0),\
\tilde\be_{35}=(0,0,0;0,1,1)/\sqrt{2},\ \tilde\be_{36}=(0,0,0;0,1,-1)/\sqrt{2}\\
\hat\be_{31}=(0,1,0;0,0,0),\ \hat\be_{32}=(0,0,1;0,0,0),
\ \hat\be_{33}=(1,0,0;0,0,0)\\
\hat\be_{34}=(0,0,0;0,1,0),\
\hat\be_{35}=(0,0,0;0,0,1),\ \hat\be_{36}=(0,0,0;1,0,0)\\
\hline
\be_{41}=(1,0,0;0,0,0),\ \be_{42}=(0,1,1;0,0,0)/\sqrt{2},
\ \be_{43}=(0,1,-1;0,0,0)/\sqrt{2}\\
\be_{44}=(0,0,0;1,0,0),\ \be_{45}=(0,0,0;0,1,1)/\sqrt{2},\
\be_{46}=(0,0,0;0,1,-1)/\sqrt{2}\\
\hline
\be_{51}=(1,0,0;0,0,0)/\sqrt{2},\ \be_{52}=(0,1,1;0,0,0)/\sqrt{2},\
\be_{53}=(0,1,-1;0,0,0)/\sqrt{2},\\
\be_{54}=(0,0,0;1,0,0),\
\be_{55}=(0,0,0;0,1,1)/\sqrt{2},\ \be_{56}=(0,0,0;0,1,-1)/\sqrt{2}\\
\end{array}
$$
\caption{Local bases at the equilibria $\xi_j$. \label{tbBases}}
\end{table}

\begin{table}[hhh]
\begin{center}
{\small
$$
\renewcommand{\arraystretch}{1.2}
\begin{array}{ccccc}
j & \xi_j & \R^6\ominus L_j & \Delta_j & \mbox{Isotypic components}\\
\hline
1 & \rz & (0,x_2,x_3;y_1,y_2,y_3) & <s_1,r\gamma^1_{\pi/2},\gamma^3_\alpha>&
(0,x_2,x_3;0,0,0),\ (0,0,0;y_1,0,0)\\
&&&&(0,0,0;0,y_2,y_2),\ (0,0,0;0,y_2,-y_2)\\
\hline
2 & \rho^2\qw & (x_1,x_2,x_3;y_1,y_2,-y_2) & <s_1,r\gamma^1_{\pi/2},\gamma^2_{\pi}>&
(x_1,0,0;0,0,0),\ (0,x_2,x_3;0,0,0)\\
&&&&(0,0,0;y_1,0,0),\ (0,0,0;0,y_2,-y_2)\\
\hline
3 & \rw & (x_1,x_2,x_3;0,y_2,y_3) & <s_1,\gamma^1_\alpha,r\gamma^2_{\pi/2}>&
(x_1,0,0;0,0,0),\ (0,x_2,x_3;0,0,0),\\
&&&& (0,0,0;0,y_2,y_3)\\
{} &\rho \rw & (x_1,x_2,x_3;y_1,y_2,0) & \rho<s_1,\gamma^1_\alpha,r\gamma^2_{\pi/2}>\rho^2&
(x_1,0,0;0,0,0),\ (0,x_2,x_3;0,0,0),\\
&&&& (0,0,0;y_1,y_2,0)\\
\hline
4 & \rho^2\qtilw & (x_1,x_2,x_3;y_1,y_2,y_2) & <rs_1,r\gamma^1_{\pi/4},\gamma^2_{\pi/2}>&
(x_1,0,0;0,0,0),\ (0,x_2,x_3;0,0,0),\\
&&&& (0,0,0;y_1,0,0),\ (0,0,0;0,y_2,y_2)\\
\hline
5 &  \rho^2\qz & (x_1,x_2,-x_2;y_1,y_2,y_3) & <s_1,r\gamma^3_{\pi/3}>&
(x_1,0,0;0,0,0),\ (0,x_2,-x_2;0,0,0),\\
&&&&(0,0,0;y_1,y_2,y_2),\ (0,0,0;0,y_2,-y_2)\\
\end{array}
$$
}
\end{center}
\caption{\label{tbIsotypic}
Isotypic decompositions of $\R^6\ominus L_j$ under $\Delta_j$.
}
\end{table}

\begin{table}[hhh]
\begin{center}
$$
\renewcommand{\arraystretch}{1.2}
\begin{array}{cccc}
\xi_i\to\xi_j & \R^6\ominus P_{ij} & \Sigma_{ij} & \mbox{Isotypic components}\\
\hline
\xi_1\to\xi_2 & (0,x_2,x_3;y_1,y_2,-y_2) & <s_1,r\gamma^1_{\pi/2},\gamma^3_{\pi}>&
(0,x_2,x_2;0,0,0),\ (0,x_2,-x_2;0,0,0)\\
&&&(0,0,0;y_1,0,0),\ (0,0,0;0,y_2,-y_2)\\
\hline
\xi_2\to\rho\xi_3 & (x_1,x_2,x_3;y_1,0,0) & <r\gamma^1_{\pi/2},r\gamma^3_{\pi/2}>&
(x_1,0,0;0,0,0),\ (0,x_2,x_3;0,0,0),\\
&&& (0,0,0;y_1,0,0)\\
\hline
\xi_3\to\xi_1 & (0,x_2,x_3;0,y_2,y_3) & <s_1,\gamma^1_{\pi},\gamma^3_{\pi/2}>&
(0,x_2,x_2;0,0,0),\ (0,x_2,-x_2;0,0,0)\\
&&&(0,0,0;0,y_2,y_2),\ (0,0,0;0,y_2,-y_2)\\
\hline
\xi_1\to\xi_4 & (0,x_2,x_3;y_1,y_2,y_2) & <r\gamma^1_{\pi/2},s_1\gamma^3_{\pi/2}>&
(0,x_2,x_3;0,0,0),\ (0,0,0;y_1,0,0),\\
&&&(0,0,0;0,y_2,y_2)\\
\hline
\xi_4\to\xi_5 & (x_1,x_2,-x_2;y_1,y_2,y_2) & <s_1r\gamma^1_{\pi},s_1r\gamma^3_{\pi}>&
(x_1,0,0;0,0,0),\ (0,x_2,-x_2;y_1,y_2,y_2)\\
\hline
\xi_5\to\rho\xi_1 & (x_1,0,0;y_1,y_2,y_3) & <r\gamma^3_{\pi/3}>&
(x_1,0,0;0,0,0),\ (0,0,0;y_1,y_2,y_3)\\
\hline
\end{array}
$$
\end{center}
\caption{\label{tbIsotypic2}
Isotypic decompositions of $\R^6\ominus P_{ij}$ under $\Sigma_{ij}$.
The plane $P_{43}$ coincides with $P_{23}$, therefore it is not listed.
}
\end{table}



\begin{table}[hhh]
{\small
\begin{center}
$$
\begin{array}{ccccc}
\mbox{Name} & \mbox{Subspace}& \mbox{Eigenvectors}& \mbox{Eigenvalues}&\mbox{Eigenvalues}\\
\hline
\rz =\xi_1& (q,0,0;0,0,0) &\be_{11}& \lambda_{11} \mbox{ (radial)}& 2A_1x^2\\
    & (0,u_2,u_3;0,0,0) &\be_{12}, \be_{13}&
   \lambda_{12}= \lambda_{13}&(A_2-A_1)x^2\\
    & (0,0,0;u,0,0) &\be_{14}&\lambda_{14}&\lambda_2+C_4x^2\\
    & (0,0,0;0,u,u) & \be_{15}&\lambda_{15}&
     \lambda_2+(C_5+C_6)x^2\\
    & (0,0,0;0,u,-u) & \be_{16}&\lambda_{16}&
     \lambda_2+(C_5-C_6)x^2\\
\hline
\rho^2\qz=\xi_5 & (q,0,0;0,0,0) &\be_{51}&\lambda_{51}& (A_2-A_1)x^2\\
    & (0,q,q;0,0,0) &\be_{52}&\lambda_{52} \mbox{ (radial)}& 2(A_1+A_2)x^2\\
    & (0,q,-q;0,0,0) & \be_{53}&\lambda_{53}& 2(A_1-A_2)x^2\\
    & (0,0,0;u_1,u_2,u_2) & \be_{54}, \be_{55}&\lambda_{54}, \lambda_{55}&\mu_1+\mu_2=2\lambda_2+(C_4+3C_5)x^2\\
    &                    &            &     & \mu_1\mu_2=(\lambda_2+(C_4+C_5)x^2)\times\\
    &                    &            &     & (\lambda_2+2C_5x^2)-2C_6^2x^4\\
    & (0,0,0;0,u_2,-u_2) & \be_{56}&\lambda_{56}&\lambda_2+(C_4+C_5)x^2\\
\hline
\rw=\xi_3 & (0,0,0;q,0,0) &\be_{34}& \lambda_{34} \mbox{ (radial)}& 2B_1y^2\\
    & (0,0,0;0,u_2,u_3) &\be_{35}, \be_{36}&\lambda_{35}, \lambda_{36}& (B_2-B_1)y^2\\
    & (u,0,0;0,0,0)   &  \be_{31}& \lambda_{31}&\lambda_1+C_1y^2\\
    & (0,u_2,u_3;0,0,0)   &  \be_{32}, \be_{33}&\lambda_{32}= \lambda_{33}&\lambda_1+C_2y^2\\
\hline
\rho\rw=\rho\xi_3 & (0,0,0;0,q,0) &\be_{34}& \lambda_{34} \mbox{ (radial)}& 2B_1y^2\\
    & (0,0,0;u_1,u_2,0) &\be_{35}, \be_{36}&\lambda_{35}, \lambda_{36}& (B_2-B_1)y^2\\
    & (0,u,0;0,0,0)   &  \be_{31}& \lambda_{31}&\lambda_1+C_1y^2\\
    & (u_1,u_2,0;0,0,0)   &  \be_{32}, \be_{33}&\lambda_{32}= \lambda_{33}&\lambda_1+C_2y^2\\
\hline
\rho^2\qw=\xi_2 & (0,0,0;u,0,0) & \be_{24}& \lambda_{24}& (B_2-B_1)y^2\\
    & (0,0,0;0,q,q) &\be_{25}& \lambda_{25} \mbox{ (radial)}& 2 (B_1+B_2)y^2\\
    & (0,0,0;0,q,-q) &\be_{26}& \lambda_{26}& 2(B_1-B_2)y^2\\
    & (q,0,0;0,0,0) & \be_{21}& \lambda_{21}&\lambda_1+(2C_2+C_3)y^2\\
    & (0,u_2,u_3;0,0,0) & \be_{22}, \be_{23}&\lambda_{22}= \lambda_{23}&\lambda_1+(C_1+C_2)y^2\\
\hline
\rho^2\qtilw=\xi_4 &  (0,0,0;0,q,q) &\be_{45}& \lambda_{45} &2(B_1-B_2)y^2\\
    & (0,0,0;0,q,-q) &\be_{46}& \lambda_{46} \mbox{ (radial)}&2(B_1+B_2)y^2\\
    & (q,0,0;0,0,0) & \be_{41}& \lambda_{41}&\lambda_1+(2C_2-C_3)y^2\\
 &(0,u_1,u_2;0,0,0)& \be_{42}, \be_{43}&\lambda_{42}=\lambda_{43}&\lambda_1+(C_1+C_2)y^2\\
    &(0,0,0;q,0,0)&\be_{44}&\lambda_{44}&(B_2-B_1)y^2
\end{array}
$$
\end{center}
}
\caption{\label{tab12w}
Eigenspaces and associated eigenvalues near equilibria in the network.
}
\end{table}

\end{document}